\documentclass[10t]{article}

\usepackage{amssymb,amsmath,amsthm,mathrsfs}
\usepackage[dvips]{graphics,epsfig}
\usepackage{verbatim}
\usepackage{colordvi,color}
\usepackage{psfrag}
\usepackage{colordvi,color,pspicture}
\usepackage{url}

\graphicspath{{Figures/}}

\usepackage[authoryear,sort&compress]{natbib}

\newcommand{\E}[1]{\mathbf{E}\,#1}

\newcommand{\Var}[1]{\mathbf{Var}\,#1}

\newcommand{\I}{\mathbf{I}}

\newcommand{\C}{\mathcal{C}}

\newcommand{\U}{\mathcal{U}}

\newcommand{\g}{\gamma}
\newcommand{\G}{\Gamma}
\newcommand{\V}{\mathcal{V}}
\newcommand{\A}{\mathcal{A}}
\newcommand{\Y}{\mathcal{Y}}
\newcommand{\y}{\mathsf{y}}
\newcommand{\X}{\mathcal{X}}
\newcommand{\M}{\mathcal{M}}

\newcommand{\mS}{\mathcal{S}}

\newcommand{\NY}{N_{\Y}}
\newcommand{\NS}{N_S}

\newcommand{\NPE}{N_{PE}}
\newcommand{\TY}{T(\Y_3)}

\newcommand{\R}{\mathbb{R}}
\newcommand{\T}{\mathcal{T}}

\newcommand{\RS}{\mathscr{R}_S}
\newcommand{\UT}{\U(T(\Y_3))}
\newcommand{\Tr}{\mathscr T_r}
\newcommand{\ve}{\varepsilon}

\DeclareMathOperator{\argmax}{argmax}
\DeclareMathOperator{\argmin}{argmin}
\DeclareMathOperator{\arginf}{arginf}

\DeclareMathOperator{\BIN}{BIN}
\DeclareMathOperator{\BER}{BER}

\theoremstyle{plain}
\newtheorem{theorem}{Theorem}[section]
\newtheorem{lemma}[theorem]{Lemma}
\newtheorem{proposition}[theorem]{Proposition}
\newtheorem{corollary}[theorem]{Corollary}

\theoremstyle{definition}

\theoremstyle{remark}
\newtheorem{conjecture}[theorem]{Conjecture}
\newtheorem{remark}[theorem]{Remark}

\begin{document}

\title{On the Distribution of the Domination Number of a New Family of Parametrized Random Digraphs$^\star$}
\author{Elvan Ceyhan\thanks{Department of Mathematics,
Ko\c{c} University, Sar{\i}yer, 34450, Istanbul, Turkey } ~~\& Carey
E. Priebe\thanks{Department of Applied Mathematics and Statistics,
The Johns Hopkins University, Baltimore, MD, 21218, USA} }
\date{\today}
\maketitle

\begin{abstract}
We derive the asymptotic distribution of the domination number of a
new family of random digraph called proximity catch digraph (PCD),
which has application to statistical testing of spatial point
patterns and to pattern recognition.
The PCD we use
is a parametrized digraph based on two sets of points on the plane,
where sample size and locations of the elements of one is held fixed,
while the sample size of the other whose elements are randomly
distributed over a region of interest goes to infinity.
PCDs are constructed based on the
relative allocation of the random set of points with respect to the
Delaunay triangulation of the other set whose size
and locations are fixed.
We introduce various auxiliary tools and concepts for the derivation of the
asymptotic distribution.
We investigate these concepts in one
Delaunay triangle on the plane, and then extend them to the multiple triangle case.
The methods are illustrated for planar data, but are
applicable in higher dimensions also.
\end{abstract}


\vspace{1.25 in}
\indent
$^\star$
This research was supported by the
Defense Advanced Research Projects Agency
as administered by the Air Force Office of Scientific Research
under contract DOD F49620-99-1-0213
and by
Office of Naval Research Grant N00014-95-1-0777.\\
\indent
$^*$Corresponding author.\\
\indent {\it E-mail address:} elceyhan@ku.edu.tr (E.~Ceyhan)

\newpage
\section{Introduction}
\label{sec:introduction}
The \emph{proximity catch digraphs} (PCDs)
are a special type of \emph{proximity graphs} which were introduced by
\cite{toussaint:1980}.
A \emph{digraph} is a directed graph with
vertices $V$ and arcs (directed edges) each of which is from one
vertex to another based on a binary relation. Then the pair $(p,q)
\in V \times V$ is an ordered pair which stands for an arc (directed
edge) from vertex $p$ to vertex $q$. For example, the \emph{nearest
neighbor (di)graph} of \cite{paterson:1992} is a proximity digraph.
 The nearest neighbor digraph has the vertex set $V$
and $(p,q)$ as an arc iff $q$ is a nearest neighbor of $p$.

Our PCDs are based on the proximity maps which are defined in a
fairly general setting. Let $(\Omega,\M)$ be a measurable space.
The \emph{proximity map} $N(\cdot)$ is defined as
$N:\Omega \rightarrow 2^\Omega$, where $2^\Omega$ is the power set of $\Omega$.
The \emph{proximity region} of $x \in \Omega$, denoted
$N(x)$, is the image of $ x \in \Omega$ under $N(\cdot)$.
The points in $N(x)$ are thought of as being ``closer" to $x \in \Omega$ than
are the points in $\Omega \setminus N(x)$.
Hence the term ``proximity" in the name \emph{proximity catch digraph}.
Proximity maps are the building blocks of the \emph{proximity graphs} of
\cite{toussaint:1980}; an extensive survey on proximity maps and graphs
is available in \cite{jaromczyk:1992}.

The \emph{proximity catch digraph} $D$ has the vertex set
$\V=\bigl\{ p_1,\ldots,p_n \bigr\}$; and the arc set $\A$ is defined
by $(p_i,p_j) \in \A$ iff $p_j \in N(p_i)$ for $i\not=j$. Notice
that the proximity catch digraph $D$ depends on the \emph{proximity}
map $N(\cdot)$ and if $p_j \in N(p_i)$, then we call $N(p_i)$ (and
hence point $p_i$) \emph{catches} $p_j$. Hence the term ``catch" in
the name \emph{proximity catch digraph}. If arcs of the form
$(p_j,p_j)$ (i.e., loops) were allowed, $D$ would have been called a
\emph{pseudodigraph} according to some authors (see, e.g.,
\cite{chartrand:1996}).

In a digraph $D=(\V,\A)$, a vertex $v \in \V$ \emph{dominates}
itself and all vertices of the form $\{u: (v,u) \in \A\}$. A
\emph{dominating set} $S_D$ for the digraph $D$ is a subset of
$\V$ such that each vertex $v \in \V$ is dominated by a vertex in
$S_D$. A \emph{minimum dominating set} $S^*_{D}$ is a dominating
set of minimum cardinality and the \emph{domination number}
$\g(D)$ is defined as $\g(D):=|S^*_{D}|$ (see, e.g., \cite{lee:1998}) where
$|\cdot|$ denotes the set cardinality functional. See
\cite{chartrand:1996} and \cite{west:2001} for more on graphs and digraphs.
If a minimum dominating set is of size one, we call it a \emph{dominating point}.

Note that for $|\V|=n>0$, $1 \le \g(D) \le n$, since $\V$ itself is always
a dominating set.

In recent years, a new classification tool based on the relative allocation of points
from various classes has been developed.
\cite{priebe:2001} introduced the \emph{class cover catch digraphs}
(CCCDs) and gave the exact and the asymptotic distribution
of the domination number of the CCCD
based on two sets, $\X_n$ and $\Y_m$, which are of size $n$ and
$m$, from classes, $\X$ and $\Y$, respectively,
and are sets of iid random variables from uniform
distribution on a compact interval in $\R$. \cite{devinney:2006},
\cite{devinney:2002a}, \cite{marchette:2003}, \cite{priebe:2003b, priebe:2003a}
applied the concept in higher dimensions and
demonstrated relatively good performance of CCCD in classification.
The methods employed involve \emph{data reduction (condensing)} by
using approximate minimum dominating sets as \emph{prototype sets}
(since finding the exact minimum dominating set is an NP-hard
problem in general
--- e.g., for CCCD in multiple dimensions --- (see \cite{devinney:2003}).
\cite{devinney:2002b} proved a SLLN result for the domination number
of CCCDs for one-dimensional data.
Although intuitively appealing and easy to extend to higher dimensions, exact
and asymptotic distribution of the domination number of the CCCDs
are not analytically tractable in $\R^2$ or higher dimensions.
As alternatives to CCCD, two new families of PCDs are introduced in
\cite{ceyhan:2003a, ceyhan:2005e} and are applied in
testing spatial point patterns (see, \cite{ceyhan:2005g, ceyhan:2006a}).
These new families are both applicable to pattern classification also.
They are designed to have better distributional and mathematical properties.
For example, the distribution of the relative density
(of arcs) is derived for one family in \cite{ceyhan:2005g} and for
the other family in \cite{ceyhan:2006a}.
In this article, we derive
the asymptotic distribution of the domination number of the latter
family called \emph{$r$-factor proportional-edge PCD}.
During the derivation process, we introduce auxiliary tools, such as,
\emph{proximity region} (which is the most crucial concept in defining
the PCD), \emph{$\G_1$-region, superset region, closest edge
extrema, asymptotically accurate distribution,} and so on.
We utilize these special regions, extrema, and asymptotic expansion of the
distribution of these extrema.
The choice of the change of
variables in the asymptotic expansion is also dependent
on the type of the extrema used and
crucial in finding the limits of the improper integrals we encounter.
Our methodology is instructive in finding the distribution of the
domination number of similar PCDs in $\R^2$ or higher dimensions.

In addition to the mathematical tractability and applicability to
testing spatial patterns and classification, this new family of PCDs
is more flexible as it allows choosing an optimal parameter for best
performance in hypothesis testing or pattern classification.

The domination number of PCDs is first investigated for data in one
Delaunay triangle (in $\R^2$) and the analysis is generalized to
data in multiple Delaunay triangles. Some trivial proofs are
omitted, shorter proofs are given in the main body of the article;
while longer proofs are deferred to the Appendix.

\section{Proximity Maps and the Associated PCDs}

\begin{figure}[ht]
\centering
 \rotatebox{-90}{ \resizebox{3. in}{!}{ \includegraphics{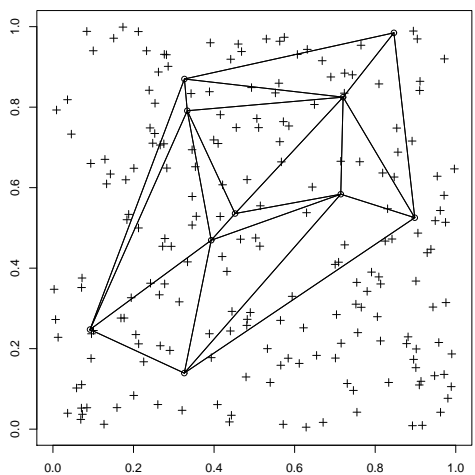}}}
 \caption{ \label{fig:deltri}
A realization of $200$ $\X$ points (crosses) and the Delaunay
triangulation based on 10 $\Y$ points (circles). }
\end{figure}

We construct the proximity regions using two data sets $\X_n$ and $\Y_m$
from two classes $\X$ and $\Y$, respectively.
Given $\Y_m \subseteq \Omega$, the {\em proximity map} $\NY(\cdot):
\Omega \rightarrow 2^{\Omega}$ associates a {\em proximity region}
$\NY(x) \subseteq \Omega$ with each point $x \in \Omega$.
The region $\NY(x)$ is defined in terms of the distance between $x$ and $\Y_m$.
More specifically, our $r$-factor proximity maps will be based on
the relative position of points from $\X_n$ with respect to the
Delaunay tessellation of $\Y_m$.
In this article, a triangle refers
to the closed region bounded by its edges.
See Figure \ref{fig:deltri} for an example with $n=200$ $\X$ points iid $\U
\bigl((0,1)\times (0,1)\bigr)$, the uniform distribution on the unit
square and the Delaunay triangulation is based on $m=10$ $\Y$
which are points also iid $\U \bigl((0,1)\times (0,1)\bigr)$.

If $\X_n=\bigl\{ X_1,\ldots,X_n \bigr\}$ is a set of $\Omega$-valued
random variables then $\NY(X_i)$ are random sets.
If $X_i$ are iid then so are the random sets $\NY(X_i)$.
We define the data-random proximity catch digraph $D$
--- associated with $\NY(\cdot)$ --- with vertex set
$\X_n=\{X_1,\cdots,X_n\}$ and arc set $\A$ by
 $$(X_i,X_j) \in \A \iff X_j \in \NY(X_i).$$
Since this relationship is not symmetric, a digraph is used rather than a graph.
The random digraph $D$ depends on the (joint)
distribution of $X_i$ and on the map $\NY(\cdot)$.
For $\X_n=\bigl\{ X_1,\cdots,X_n \bigr\}$, a set of iid random variables
from $F$, the domination number of the associated data-random
proximity catch digraph based on the proximity map $N(\cdot)$,
denoted $\g(\X_n,N)$, is the minimum number of point(s) that dominate
all points in $\X_n$.

The random variable $\g(\X_n,N)$ depends explicitly on $\X_n$ and $N(\cdot)$
and implicitly on $F$ .
Furthermore, in general, the
distribution, hence the expectation $\E[\g(\X_n,N)]$, depends on
$n$, $F$, and $N$; $ 1 \leq \E[\g(\X_n,N)] \le n.$
In general, the variance of $\g(\X_n,N)$ satisfies,
$1 \le \Var[\g(\X_n,N)] \le n^2/4$.

For example, the CCCD of \cite{priebe:2001} can be viewed as an
example of PCDs and is briefly discussed in the next section.
We use many of the properties of CCCD in $\R$ as guidelines in
defining PCDs in higher dimensions.

\subsection{Spherical Proximity Maps}
 Let $\Y_m=\left \{\y_1,\ldots,\y_m \right\} \subset \R$.
Then the proximity map associated with CCCD is defined as the open
ball $\NS(x):=B(x,r(x))$ for all $x \in \R$, where $r(x):=\min_{\y
\in \Y_m}d(x,\y)$ (see \cite{priebe:2001}) with $d(x,y)$ being the
Euclidean distance between $x$ and $y$.
That is, there is an arc
from $X_i$ to $X_j$ iff there exists an open ball centered at $X_i$
which is ``pure" (or contains no elements) of $\Y_m$ in its
interior, and simultaneously contains (or ``catches") point $X_j$.
We consider the closed ball, $\overline{B}(x,r(x))$ for $\NS(x)$ in this article.
Then for $x \in \Y_m$, we have $\NS(x)=\{x\}$.
Notice that a
ball is a sphere in higher dimensions, hence the notation $N_S$.
Furthermore, dependence on $\Y_m$ is through $r(x)$.
Note that  in $\R$ this
proximity map is based on the intervals
$I_j=\left(\y_{(j-1):m},\y_{j:m} \right)$ for $j=0,\ldots,m+1$ with
$\y_{0:m}=-\infty$ and $\y_{(m+1):m}=\infty$, where $\y_{j:m}$ is the
$j^{th}$ order statistic in $\Y_m$. This interval partitioning can
be viewed as the Delaunay tessellation of $\R$ based on $\Y_m$. So
\emph{in higher dimensions, we use the Delaunay triangulation based
on $\Y_m$ to partition the support}.

A natural extension of the proximity region $\NS(x)$ to $\R^d$ with
$d>1$ is obtained as $N_S(x):=B(x,r(x))$ where
$r(x):=\min_{\y \in \Y_m} d(x,\y)$ which
is called the \emph{spherical proximity map}.
The spherical proximity map $N_S(x)$ is well-defined for all $x \in \R^d$ provided
that $\Y_m \not= \emptyset$.
Extensions to $\R^2$ and higher
dimensions with the spherical proximity map --- with applications in
classification --- are investigated by  \cite{devinney:2006},
\cite{devinney:2002a}, \cite{marchette:2003}, \cite{priebe:2003b, priebe:2003a}.
However, finding the minimum dominating set
of CCCD (i.e., the PCD associated with $N_S(\cdot)$) is an NP-hard
problem and the distribution of the domination number is not
analytically tractable for $d>1$.
This drawback has motivated us to define
new types of proximity maps.
\cite{ceyhan:2005e} introduced $r$-factor proportional-edge PCD,
where the distribution
of the domination number of $r$-factor PCD with $r=3/2$ is used in
testing spatial patterns of segregation or association.
\cite{ceyhan:2006a} computed the asymptotic distribution of the
relative density of the $r$-factor PCD and used it for the same purpose.
\cite{ceyhan:2003a} introduced the central similarity
proximity maps and the associated PCDs, and \cite{ceyhan:2005g}
computed the asymptotic distribution of the relative density of the
parametrized version of the central similarity PCDs and applied the
method to testing spatial patterns. An extensive treatment of the
PCDs based on Delaunay tessellations is available in \cite{ceyhan:2004d}.

The following property (which is referred to as
\textbf{Property (\ref{prop:P6})}) of CCCDs in $\R$ plays an important role in
defining proximity maps in higher dimensions.
\begin{equation}
\label{prop:P6}
 \text{ \textbf{Property (1)} For $x \in I_j$, $\NS(x)$ is a proper subset
of $I_j$ for almost all $x \in I_j$.}
\end{equation}
In fact, \textbf{Property (\ref{prop:P6})} holds for all $x \in
I_j\setminus \{(\y_{(j-1):m}+\y_{j:m})/2\}$ for CCCDs in $\R$.
For $x \in I_j$, $\NS(x)=I_j$ iff
$x=\left(\y_{(j-1):m}+\y_{j:m} \right)/2$.
We define an associated region for such points in the general context.
The \emph{superset region} for any proximity map $N(\cdot)$ in
$\Omega$ is defined to be
$$\RS(N):=\bigl\{ x \in \Omega: N(x) =\Omega \bigr\}.$$
For example, for $\Omega=I_j \subsetneq \R$,
$\RS(\NS):=\{x \in I_j: \NS(x) = I_j\}=
\left\{ \left(\y_{(j-1):m}+\y_{j:m} \right)/2 \right\}$
and for $\Omega=\T_j \subsetneq \R^d$,
$\RS(\NS):=\{x \in \T_j: \NS(x) = \T_j\},$
where $\T_j$ is the $j^{th}$ Delaunay cell in the Delaunay
tessellation. Note that for $x \in I_j$, $\lambda(\NS(x)) \le
\lambda(I_j)$ and $\lambda(\NS(x)) = \lambda(I_j)$ iff $x \in
\RS(\NS)$ where $\lambda(\cdot)$ is the Lebesgue measure on $\R$. So
the proximity region of a point in $\RS(\NS)$ has the largest Lebesgue measure.
Note also that given $\Y_m$, $\RS(\NS)$ is
not a random set, but $\I(X\in \RS(\NS))$ is a random variable,
where $\I(\cdot)$ stands for the indicator function.
\textbf{Property (\ref{prop:P6})} also implies that $\RS(\NS)$ has
zero $\R$-Lebesgue measure.

Furthermore, given a set $B$ of size $n$ in $[\y_{1:m},\y_{m:m}] \setminus \Y_m$,
the number of disconnected components in the PCD based on $N_S(\cdot)$
is at least the cardinality of the set
$\{j \in \{1,2,\ldots,m\}: B \cap I_j \not=\emptyset \}$,
which is the set of indices of the intervals that contain
some point(s) from $B$.

Since the distribution of the domination number of spherical PCD
(or CCCD) is tractable in $\R$, but not in $\R^d$ with $d>1$, we try
to mimic its properties in $\R$ while defining new PCDs in higher
dimensions.

\section{The $r$-Factor Proportional-Edge Proximity Maps}
\label{sec:r-factor} First, we describe the construction of the
$r$-factor proximity maps and regions, then state some of its basic
properties and introduce some auxiliary tools.

\subsection{Construction of the Proximity Map}
Let $\Y_m=\left \{\y_1,\ldots,\y_m \right\}$ be $m$ points in
general position in $\R^d$ and $\T_j$ be the $j^{th}$ Delaunay cell
for $j=1,\ldots,J_m$, where $J_m$ is the number of Delaunay cells.
Let $\X_n$ be a set of iid random variables from distribution $F$ in
$\R^d$ with support $\mS(F) \subseteq \C_H(\Y_m)$.

In particular, for illustrative purposes, we focus on $\R^2$ where
a Delaunay tessellation is a triangulation, provided that no more
than three points in $\Y_m$ are cocircular (i.e., lie in the same circle).
Furthermore, for simplicity,
let $\Y_3=\{\y_1,\y_2,\y_3\}$ be three non-collinear points
in $\R^2$ and $\TY=T(\y_1,\y_2,\y_3)$ be the triangle
with vertices $\Y_3$.
Let $\X_n$ be a set of iid random variables from $F$ with
support $\mS(F) \subseteq \TY$.
If $F=\U(\TY)$, a composition of translation,
rotation, reflections, and scaling
will take any given triangle $\TY$
to the basic triangle $T_b=T((0,0),(1,0),(c_1,c_2))$
with $0 < c_1 \le 1/2$, $c_2>0$,
and $(1-c_1)^2+c_2^2 \le 1$, preserving uniformity.
That is, if $X \sim \U(\TY)$ is transformed in the same manner to,
say $X'$, then we have $X' \sim \U(T_b)$.

\begin{figure} [ht]
    \centering
   \scalebox{.35}{\input{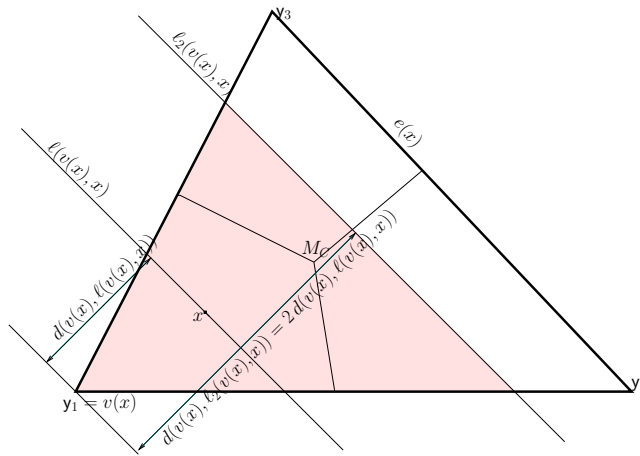}}
   \caption{Construction of $r$-factor proximity region,
$\NPE^{r=2}(x)$ (shaded region).} \label{fig:ProxMapDef1}
\end{figure}

\begin{figure}
\begin{center}
\psfrag{A}{\scriptsize{$\y_1$}}
 \psfrag{B}{\scriptsize{$\y_2$}}
\psfrag{C}{\scriptsize{$\y_3$}}
 \psfrag{CM}{\scriptsize{$M_C$}}
\psfrag{IC}{\scriptsize{$M_{I}$}}
 \psfrag{x}{}
 \psfrag{P1}{\scriptsize{$M_1$}}
 \psfrag{P2}{\scriptsize{$M_2$}}
 \psfrag{P3}{\scriptsize{$M_3$}}
\psfrag{R(A)}{\scriptsize{$R_{M_C}(A)$} }
\psfrag{R(B)}{\scriptsize{$R_{M_C}(B)$} }
\psfrag{R(C)}{\scriptsize{$R_{M_C}(C)$} }
 \epsfig{figure=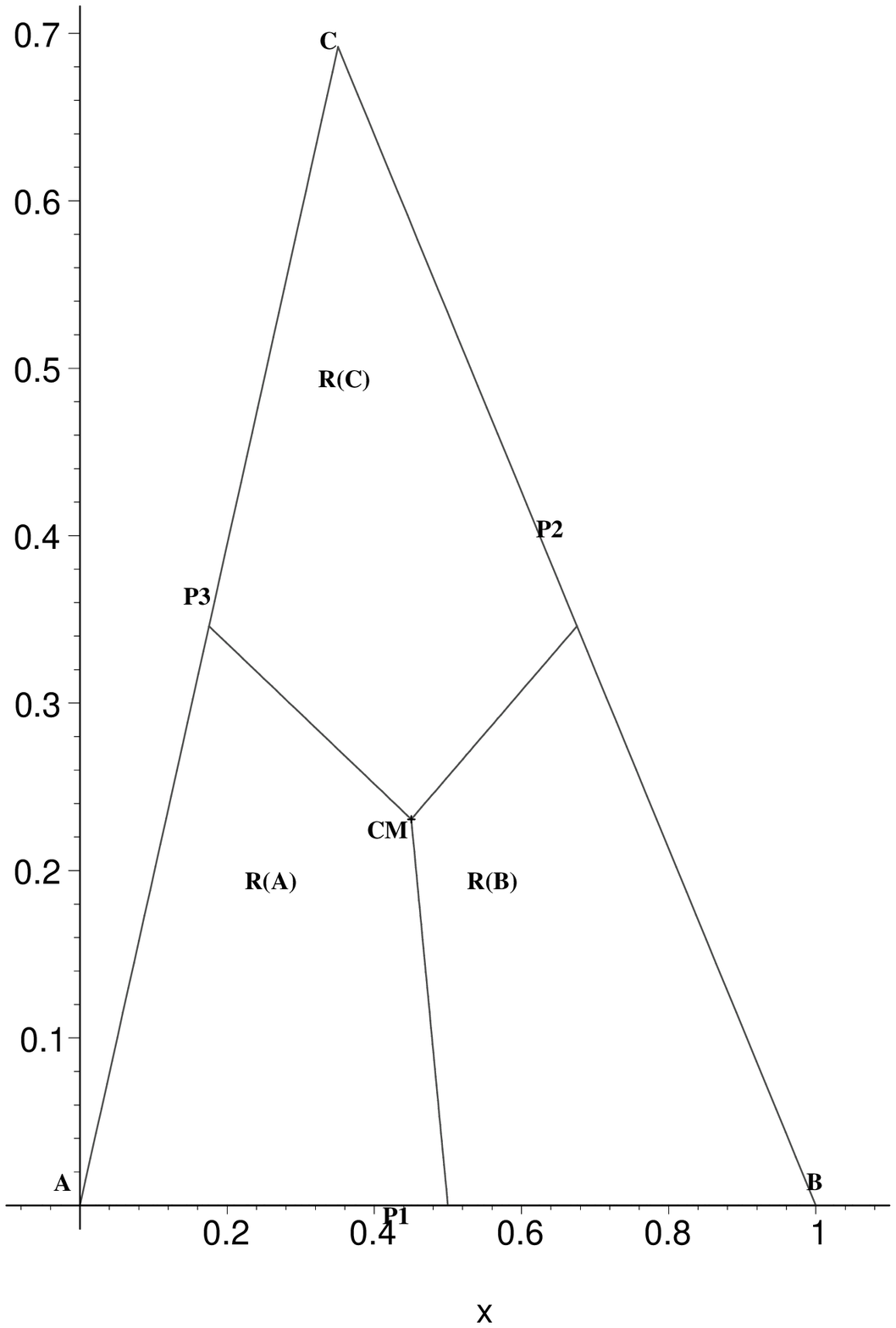, height=140pt,width=200pt}
  \psfrag{1}{\scriptsize{$\y_2$}}
\psfrag{R(A)}{\scriptsize{$R_{M_I}(A)$} }
\psfrag{R(B)}{\scriptsize{$R_{M_I}(B)$} }
\psfrag{R(C)}{\scriptsize{$R_{M_I}(C)$} }
 \epsfig{figure=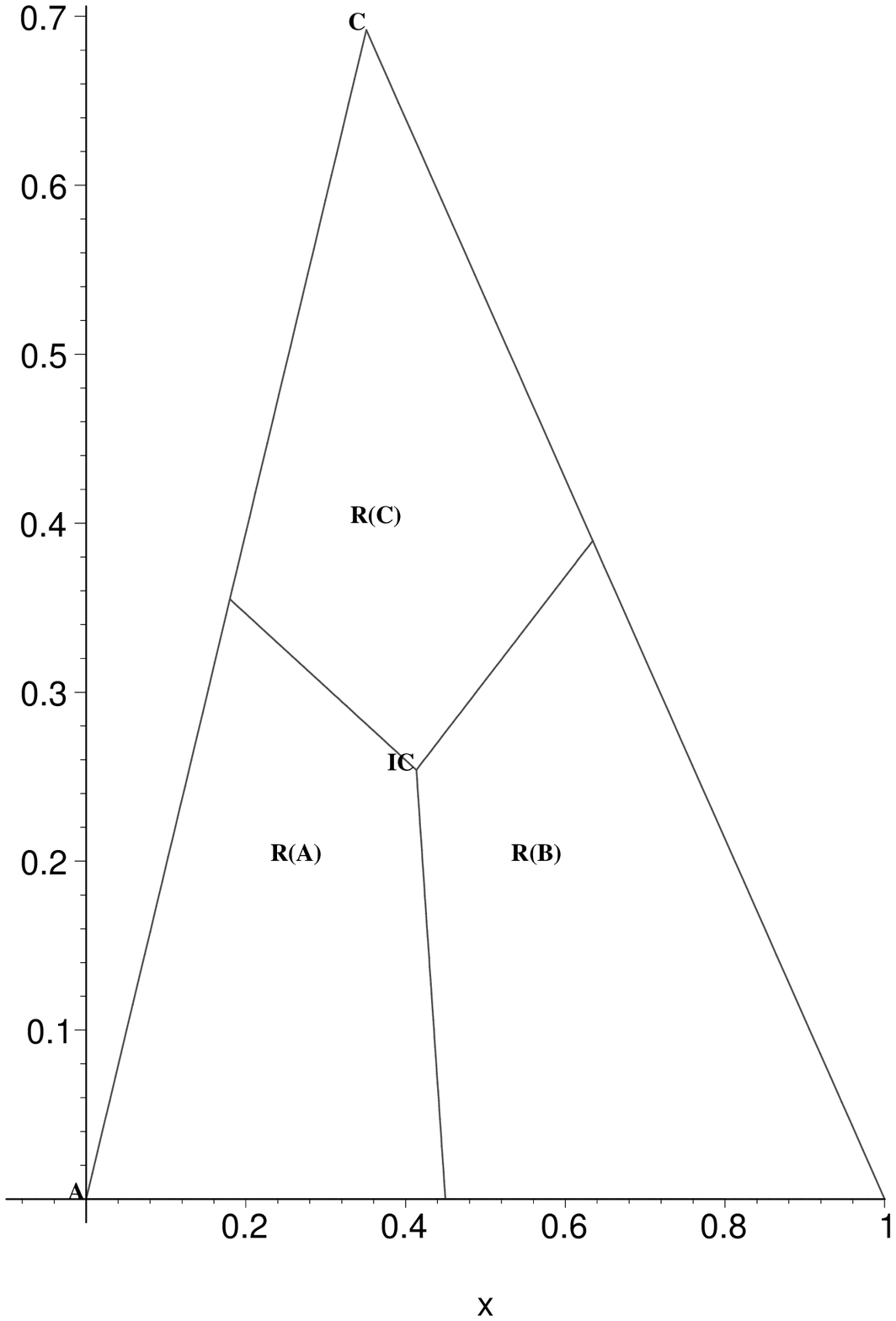, height=140pt,width=200pt}
\end{center}
\caption{
 \label{fig:vert-regions}
The vertex regions constructed with center of mass $M_C$ (left) and
incenter $M_I$ (right) using the line segments on the line joining
$M$ to $\Y_3$. }
\end{figure}

For $r \in [1,\infty]$, define $\NPE^r(\cdot,M):=N(\cdot,M;r,\Y_3)$
to be the \emph{$r$-factor proportional-edge proximity map} with
$M$-vertex regions as follows (see also Figure \ref{fig:ProxMapDef1}
with $M=M_C$ and $r=2$).
For $x \in \TY \setminus \Y_3$, let $v(x)
\in \Y_3$ be the vertex whose region contains $x$; i.e., $x \in
R_M(v(x))$. In this article \emph{$M$-vertex regions} are
constructed by the lines joining any point $M \in \R^2 \setminus
\Y_3$ to a point on each of the edges of $\TY$. Preferably, $M$ is
selected to be in the interior of the triangle $\TY^o$. For such an
$M$, the corresponding vertex regions can be defined using the line
segment joining $M$ to $e_j$, which lies on the line joining $\y_j$
to $M$; e.g. see Figure \ref{fig:vert-regions} (left) for vertex
regions based on center of mass $M_C$, and (right) incenter $M_I$.
With $M_C$, the lines joining $M$ and $\Y_3$ are the \emph{median
lines}, that cross edges at $M_j$ for $j=1,2,3$. $M$-vertex regions,
among many possibilities, can also be defined by the orthogonal
projections from $M$ to the edges. See \cite{ceyhan:2004d} for a
more general definition.
The vertex regions in Figure \ref{fig:ProxMapDef1} are
center of mass vertex regions or $CM$-vertex regions.
If $x$ falls on the boundary of two $M$-vertex regions,
we assign $v(x)$ arbitrarily.
Let $e(x)$ be the edge of $\TY$ opposite of $v(x)$.
Let $\ell(v(x),x)$ be the line parallel to $e(x)$ through
$x$. Let $d(v(x),\ell(v(x),x))$ be the Euclidean (perpendicular)
distance from $v(x)$ to $\ell(v(x),x)$. For $r \in [1,\infty)$, let
$\ell_r(v(x),x)$ be the line parallel to $e(x)$ such that
$$
d(v(x),\ell_r(v(x),x)) = r\,d(v(x),\ell(v(x),x))\\
\text{ and }\\
d(\ell(v(x),x),\ell_r(v(x),x)) < d(v(x),\ell_r(v(x),x)).
$$
Let $T_r(x)$ be the triangle similar to and with the same
orientation as $\TY$ having $v(x)$ as a vertex and $\ell_r(v(x),x)$
as the opposite edge. Then the {\emph r-factor proportional-edge
proximity region} $\NPE^r(x,M)$ is defined to be $T_r(x) \cap \TY$.
Notice that $\ell(v(x),x)$ divides the edges of $T_r(x)$ (other than
the one lies on $\ell_r(v(x),x)$) proportionally with the factor
$r$. Hence the name \emph{$r$-factor proportional edge proximity
region}.

Notice that $r \ge 1$ implies $x \in \NPE^r(x,M)$ for all $x \in \TY$.
Furthermore,
$\lim_{r \rightarrow \infty} \NPE^r(x,M) = \TY$ for all $x \in \TY \setminus \Y_3$,
so we define $\NPE^{\infty}(x,M) = \TY$ for all such $x$.
For $x \in \Y_3$, we define $\NPE^r(x,M) = \{x\}$ for all $r \in [1,\infty]$.

\begin{figure}[ht]
\centering
 \rotatebox{-90}{ \resizebox{2.1 in}{!}{\includegraphics{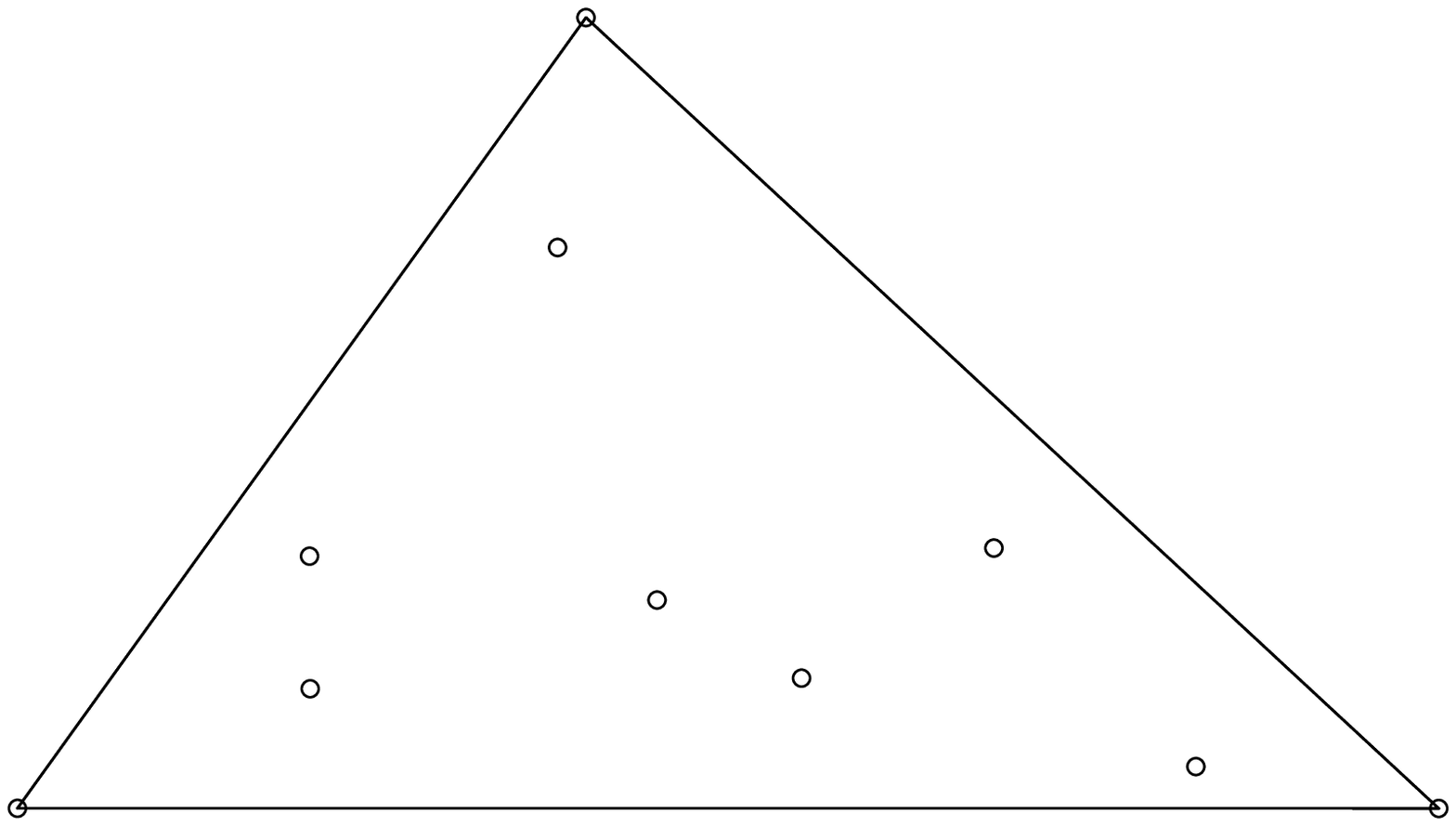}}}
 \rotatebox{-90}{ \resizebox{2.1 in}{!}{\includegraphics{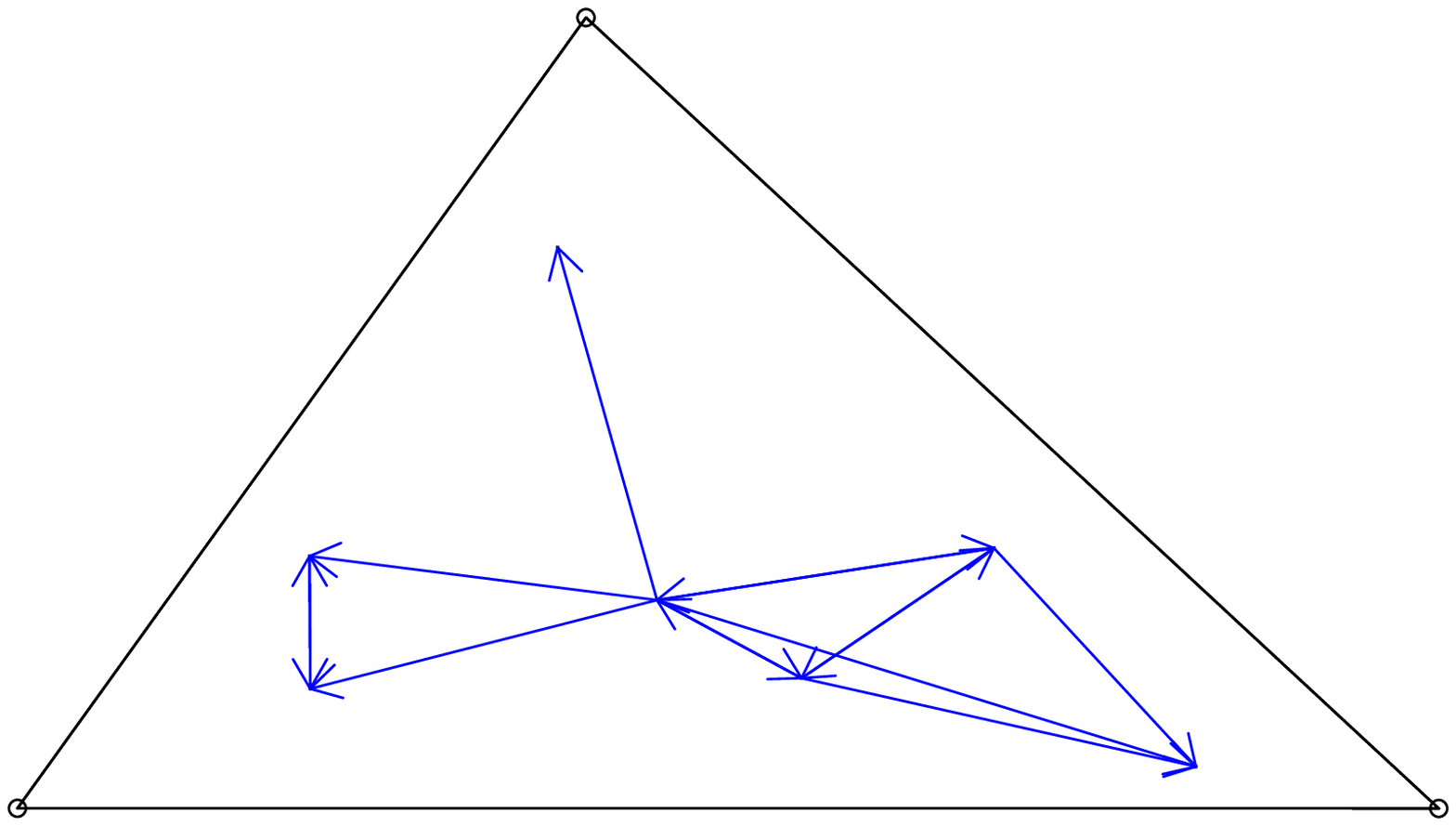}}}
 \caption{
 \label{fig:one-tri-arcs}
 A realization of 7 $\X$ points generated iid $\U(\TY)$ (left)
 and the corresponding arcs of $r$-factor proportional edge PCD with $r=3/2$ and $M=M_C$. }
\end{figure}

Hence, $r$-factor proportional edge PCD has vertices $\X_n$ and arcs
$(x_i,x_j)$ iff $x_j \in \NPE^r(x_i,M)$. See Figure
\ref{fig:one-tri-arcs} for a realization of $\X_n$ with $n=7$ and
$m=3$. The number of arcs is 12 and $\g_n(r=2,M_C)=1$.  By
construction, note that as $x$ gets closer to $M$ (or equivalently
further away from the vertices in vertex regions), $\NPE^r(x,M)$
increases in area, hence it is more likely for the outdegree of $x$
to increase.  So if more $\X$ points are around the center $M$, then
it is more likely for $\g_n$ to decrease, on the other hand, if more
$\X$ points are around the vertices $\Y_3$, then the regions get
smaller, hence it is more likely for the outdegree for such points
to be smaller, thereby implying $\g_n$ to increase.  This
probabilistic behaviour is utilized in \cite{ceyhan:2005e} for
testing spatial patterns.

Note also that, $\NPE^r(x,M)$ is a \emph{homothetic transformation
(enlargement)} with $r \ge 1$ applied on the region
$\NPE^{r=1}(x,M)$. Furthermore, this transformation is also an
\emph{affine similarity transformation}.

\subsection{Some Basic Properties and Auxiliary Concepts}

First, notice that $X_i \stackrel{iid}{\sim} F$, with the additional
assumption that the non-degenerate two-dimensional probability
density function $f$ exists with support $\mS(F) \subseteq \TY$,
imply that the special case in the construction of $\NPE^r$ ---
$X$ falls on the boundary of two vertex regions --- occurs with
probability zero. Note that for such an $F$, $\NPE^r(X)$ is a
triangle a.s.

The similarity ratio of $\NPE^r(x,M)$ to $\TY$ is given by
$\frac{\min \Bigl(d\bigl(v(x),\,e(x)\bigr),r\,d\bigl(v(x),\,
\ell(v(x),x)\bigr)\Bigr)} {d(v(x),\,e(x))},$
 that is, $\NPE^r(x,M)$ is similar to $\TY$ with
the above ratio. \textbf{Property (\ref{prop:P6})} holds depending
on the pair $M$ and $r$. That is, there exists an $r_0$ and a
corresponding point $M(r_0) \in \TY^o$ so that $\NPE^{r_0}(x,M)$
satisfies \textbf{Property (\ref{prop:P6})} for all $r \le r_0$,
but fails to satisfy it otherwise. \textbf{Property (\ref{prop:P6})}
fails for all $M$ when $r=\infty$. With $CM$-vertex regions, for all
$r \in [1,\infty]$, the area $A\left(\NPE^r(x,M_C) \right)$ is a
continuous function of $d(\ell_r(v(x),x),v(x))$ which is a
continuous function of $d(\ell(v(x),x),v(x))$ which is a continuous
function of $x$.



Note that if $x$ is close enough to $M$,
we might have $\NPE^r(x,M) = \TY$ for $r=\sqrt{2}$ also.

\begin{figure}
\begin{center}
\psfrag{A}{\scriptsize{$\y_1$}}
 \psfrag{B}{\scriptsize{$\y_2$}}
\psfrag{C}{\scriptsize{$\y_3$}}
 \psfrag{CC}{\scriptsize{$M_{CC}$}}
\psfrag{IC}{\scriptsize{$M_{I}$}}
 \psfrag{x}{\scriptsize{$x$}}
\psfrag{D}{}
 \psfrag{E}{}
 \psfrag{F}{}
\psfrag{S}{}
 \psfrag{Ac}{\scriptsize{$q_1(r,x)$} }
  \psfrag{Ab}{}
  \psfrag{Ba}{}
\psfrag{Bc}{\scriptsize{$q_2(r,x)$} }
 \psfrag{Ca}{}
\psfrag{Cb}{\scriptsize{$q_3(r,x)$} }
 \epsfig{figure=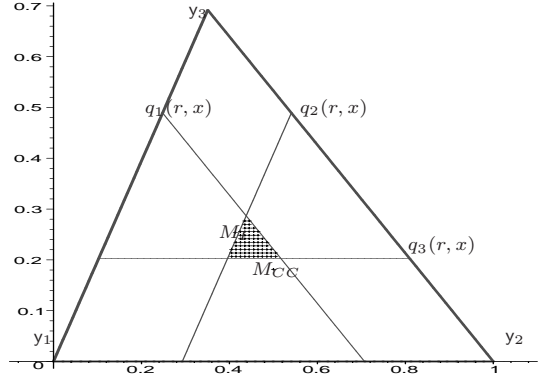, height=140pt , width=200pt}
\end{center}
\caption{The triangle $\Tr$ with $r=\sqrt{2}$ (the hatched region).
} \label{fig:Tr-RS-NDA-CC}
\end{figure}

In $\TY$, drawing the lines $q_j(r,x)$ such that
$d(\y_j,e_j)=r\,d(\y_j,q_j(r,x))$ for $j\in \{1,2,3\}$  yields a
triangle, denoted $\Tr$, for $r<3/2$ . See Figure
\ref{fig:Tr-RS-NDA-CC} for $\Tr$ with $r=\sqrt{2}$.

\begin{figure}
\begin{center}
\psfrag{A}{\scriptsize{$\y_1$}}
 \psfrag{B}{\scriptsize{$\y_2$}}
\psfrag{C}{\scriptsize{$\y_3$}}
 \psfrag{CC}{\scriptsize{$M_{CC}$}}
 \psfrag{x}{}
 \psfrag{y}{}
  \psfrag{M1}{}
 \psfrag{M2}{}
 \psfrag{M3}{}
 \epsfig{figure=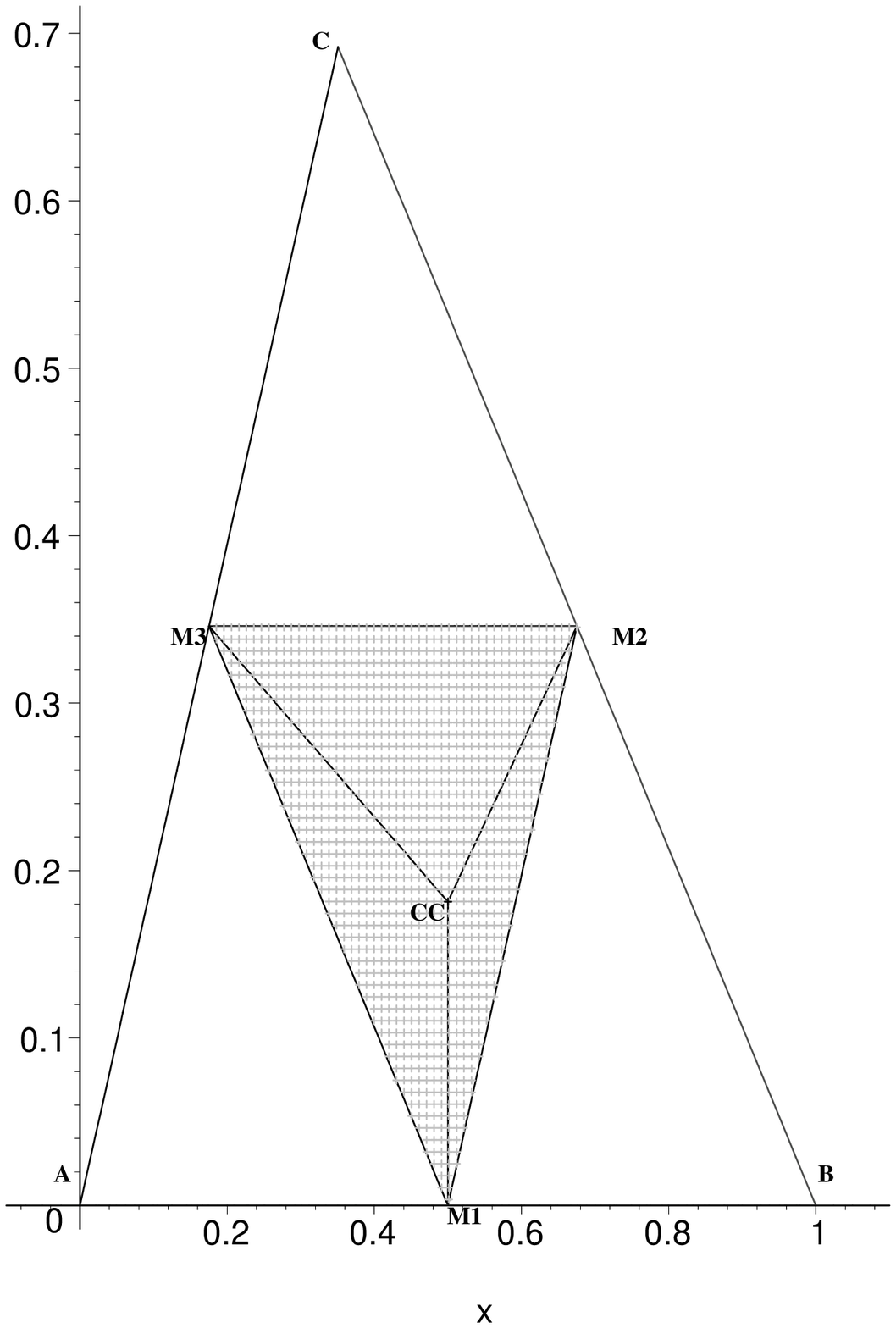, height=140pt,width=200pt}
 \epsfig{figure=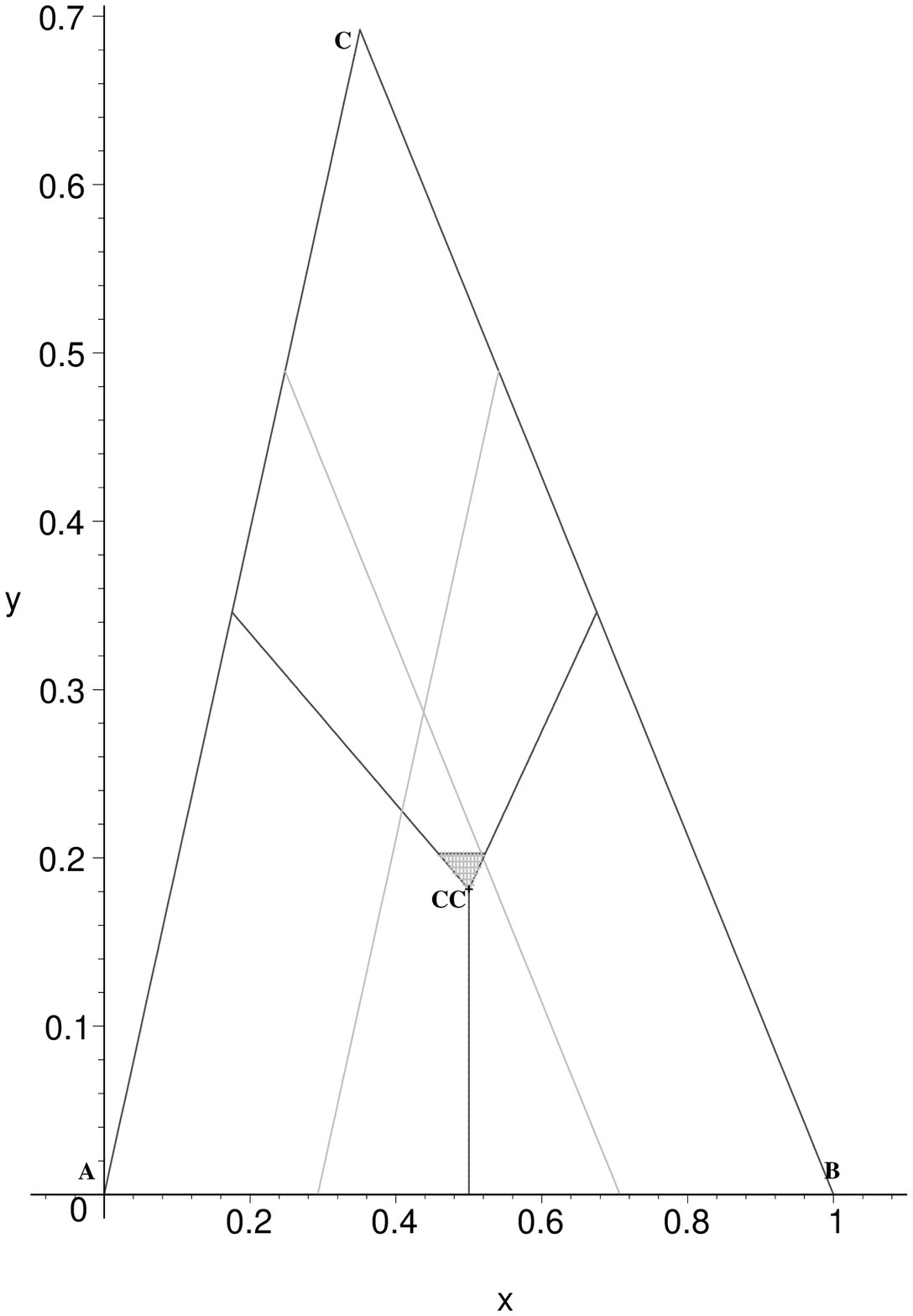, height=140pt,width=200pt}
\end{center}
\caption{
 \label{fig:SS-regions}
The superset regions (the shaded regions) constructed with
circumcenter $M_{CC}$ with $r=\sqrt{2}$ (left) and $r=2$ (right)
with vertex regions constructed with orthogonal projections to the
edges. }
\end{figure}

The functional form of $\Tr$ in the basic triangle $T_b$ is given by
{\small
\begin{align}
\label{eqn:T^r-def}
& \Tr =T(t_1(r),t_2(r),t_3(r))= \left \{(x,y) \in T_b:
y \ge \frac{c_2\,(r-1)}{r};\; y \le \frac{c_2\,(1-r\,x)}{r\,(1-c_1)};\;
y \le \frac{c_2\,(r\,(x-1)+1)}{r\,c_1} \right\}\\
&=T\Biggl( \left(\frac{(r-1)\,(1+c_1)}{r},\frac{c_2\,(r-1)}{r} \right),
\left(\frac{2-r+c_1\,(r-1)}{r},\frac{c_2\,(r-1)}{r} \right),
\left(\frac{c_1\,(2-r)+r-1}{r},\frac{c_2\,(r-2)}{r} \right) \Biggr) \nonumber.
\end{align}
}
There is a crucial difference between the triangles $\Tr$ and $T(M_1,M_2,M_3)$.
More specifically $T(M_1,M_2,M_3) \subseteq \RS(r,M)$ for all $M$ and $r \ge 2$, but
$(\Tr)^o$ and $\RS(r,M)$ are disjoint for all $M$ and $r$. So if $M
\in (\Tr)^o$, then $\RS(r,M)=\emptyset$; if $M \in
\partial(\Tr)$, then $\RS(r,M)=\{M\}$; and if $M \not\in \Tr$,
then $\RS(r,M)$ has positive area. Thus $\NPE^r(\cdot,M)$ fails to
satisfy \textbf{Property (\ref{prop:P6})} if $M \not\in \Tr$. See
Figure \ref{fig:SS-regions} for two examples of superset regions
with $M$ that corresponds to circumcenter $M_{CC}$ in this triangle
and the vertex regions are constructed using orthogonal projections.
For $r=2$, note that $\Tr=\emptyset$ and the superset region is
$T(M_1,M_2,M_3)$ (see Figure \ref{fig:SS-regions} (left)), while for
$r=\sqrt{2}$, $\Tr^o$ and $\RS(r=\sqrt{2},M)^o$ are disjoint (see
Figure \ref{fig:SS-regions} (right))

The triangle $\Tr$ given in Equation (\ref{eqn:T^r-def}) and the
superset region $\RS(r,M)$ play a crucial role in computing the
distribution of the domination number of the $r$-factor PCD.

\subsection{Main Result}
Next, we present the main result of this article. Let
$\g_n(r,M):=\g\left(\X_n, \NPE^r,M \right)$ be the domination number
of the PCD based on $\NPE^r$ with $\X_n$, a set of iid random
variables from $\UT$, with $M$-vertex regions.

The domination number $\g_n(r,M)$ of the PCD has the following
asymptotic distribution. As $n \rightarrow \infty$,
\begin{equation}
\label{eqn:asymptotic-NYr}
 \g_n(r,M) \sim
\left\lbrace \begin{array}{ll}
       2+\BER(1-p_r),           & \text{for $r \in [1,3/2]$ and $M \in \{t_1(r),t_2(r),t_3(r)\}$,}\\
       1,           & \text{for $r>3/2$,}\\
       3,           & \text{for $r \in [1,3/2)$ and $M \in \Tr\setminus \{t_1(r),t_2(r),t_3(r)\}$,}\\
\end{array} \right.
\end{equation}
where $\BER(p)$ stands for Bernoulli distribution with probability of
success $p$, $\Tr$ and $t_j(r)$ are defined in Equation
(\ref{eqn:T^r-def}), and for $r \in [1,3/2)$ and $M \in
\{t_1(r),t_2(r),t_3(r)\}$,
\begin{equation}
\label{eqn:p_r-form} p_r=\int_0^{\infty}\int_0^{\infty}\frac
{64\,r^2}{9\,(r-1)^2}\,w_1\,w_3\,\exp\left(\frac{4\,r}
{3\,(r-1)}\,(w_1^2+w_3^2+2\,r\,(r-1)\,w_1\,w_3)\right)\,dw_3w_1;
\end{equation}
for example for $r=3/2$ and $M=M_C$, $p_r \approx .7413$.

In Equation (\ref{eqn:asymptotic-NYr}), the first line is referred
as the non-degenerate case, the second and third lines are referred
as degenerate cases with a.s. limits 1 and 3, respectively.

\begin{figure}
\begin{center}
\psfrag{r}{\scriptsize{$r$}} \psfrag{pr}{\scriptsize{$p_r$}}
\epsfig{figure=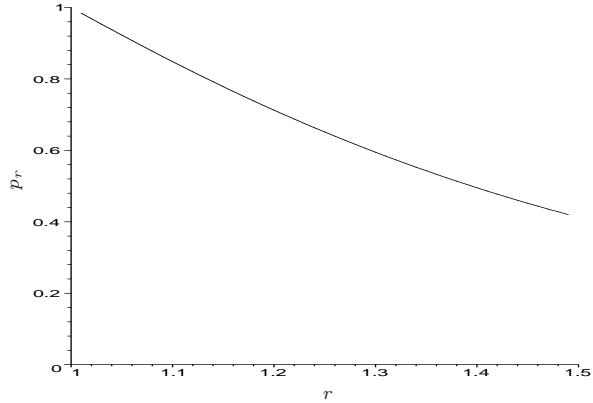, height=150pt, width=220pt}
\end{center}
\caption{ Plotted is the probability $p_r=\lim_{n\rightarrow
\infty}P\left( \g_n(r,M)=2 \right)$ given in Equation
(\ref{eqn:p_r-form}) as a function of $r$ for $r \in [1,3/2)$ and $M
\in \{t_1(r),t_2(r),t_3(r)\}$.} \label{fig:pglt2ofr}
\end{figure}

In the following sections, we define a region associated with $\g=1$
case in general. Then we give finite sample and asymptotic upper
bounds for $\g_n(r,M)$. Then we derive the asymptotic distribution
of $\g_n(r,M)$.

\section{The $\G_1$-Regions for $\NPE^r$}
\label{sec:Gamma1-NYr}
First, we define $\G_1$-regions in general,
and describe the
construction of $\G_1$-region of $\NPE^r$ for one point and
multiple point data sets, and provide some results concerning $\G_1$-regions.

\subsection{Definition of $\G_1$-Regions}
Let $(\Omega,\mathcal{M})$ be a measurable space and consider the
proximity map $N:\Omega \rightarrow 2^{\Omega}$. For any set $B
\subseteq \Omega$, the \emph{$\G_1$-region} of $B$ associated with
$N(\cdot)$, is defined to be the region $\G^N_1(B):=\{z \in \Omega:
B \subseteq  N(z)\}$. For $x \in \Omega$, we denote
$\G^N_1\bigl(\{x\}\bigr)$ as $\G^N_1(x)$.

If $\X_n=\bigl\{ X_1,X_2,\cdots,X_n \bigr\}$ is a set
of $\Omega$-valued random variables,
then $\G^N_1(X_i)$, $i=1,\cdots,n$, and $\G^N_1(\X_n)$ are random sets.
If the $X_i$ are iid,
then so are the random sets $\G^N_1(X_i)$.

 Note that $\g(\X_n,N)=1$ iff $\X_n \cap \G^N_1(\X_n)\not= \emptyset$.
Hence the name \emph{$\G_1$-region}.

It is trivial to see the following.
\begin{proposition}
\label{prop:RSsubsetG1}
For any proximity map $N$ and set $B \subseteq \Omega$, $\RS(N) \subseteq \G^N_1(B)$.
\end{proposition}

\begin{lemma}
\label{lem:gamma1-intersect} For any proximity map $N$ and $B
\subseteq \Omega$, $\G^N_1(B)= \cap_{x \in B}\G^N_1(x)$.
\end{lemma}
{\bfseries Proof:}  Given a particular type of proximity map $N$ and
subset $B \subseteq \Omega$, $y \in \G^N_1(B)$ iff $B \subseteq
N(y)$ iff $x \in N(y)$ for all $x \in B$ iff $y \in \G^N_1(x)$ for
all $x \in B$ iff $y \in \cap_{x \in B}\G^N_1(x)$. Hence the result
follows. $\blacksquare$

A problem of interest is finding, if possible, a (proper) subset of
$B$, say $G \subsetneq B$, such that $\G^N_1(B)=\cap_{x \in
G}\G^N_1(x)$. This implies that only the points in $G$ will be
\emph{active} in determining $\G^N_1(B)$.



For example, in $\R$ with $\Y_2=\{0,1\}$, and $\X_n$ a set of iid random variables of size $n>1$ from
$F$ in $(0,1)$,
 $\G^{N_S}_1(\X_n)=\Bigl(X_{n:n}/2,(1+X_{1:n})/2 \Bigr)$.
So the extrema (minimum and maximum) of the set
$\X_n$ are sufficient to determine the $\G_1$-region;
i.e., $G=\{X_{1:n},X_{n:n} \}$
for $\X_n$ a set of iid random variables from a continuous distribution on $(0,1)$.
Unfortunately, in the multi-dimensional case, there is no natural ordering that
yields natural extrema such as minimum or maximum.

\subsection{Construction of $\G_1$-Region of a Point for $\NPE^r$}
\begin{figure} [ht]
    \centering
    \scalebox{.35}{\input{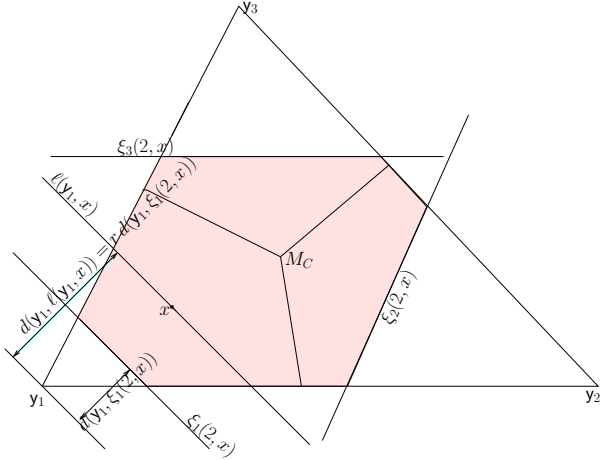}}
    \caption{Construction of the $\G_1$-region, $\G^{r=2}_1(x,M_C)$ (shaded region). }
\label{fig:ProxMapDef2}
\end{figure}

For $\NPE^r(\cdot,M)$, the $\G_1$-region,
denoted as $\G^r_1(\cdot,M):=\G^{\NPE^r}_1(\cdot,M)$, is constructed as follows;
see also Figure \ref{fig:ProxMapDef2}.
Let $\xi_j(r,x)$ be the line
parallel to $e_j$ such that $\xi_j(r,x)\cap \TY \not=\emptyset$ and
$r\,d(\y_j,\xi_j(r,x))=d(\y_j,\ell(\y_j,x))$ for $j\in \{1,2,3\}$.
Then
$$\G^r_1(x,M)=\cup_{j=1}^3 \bigl[ \G^r_1(x,M)\cap R_M(\y_j)\bigr]$$
where
$\G^r_1(x,M)\cap R_M(\y_j)=\{z \in R_M(\y_j):
d(\y_j,\ell(\y_j,z)) \ge d(\y_j,\xi_j(r,x)\} \text{ for } j\in
\{1,2,3\}.$

Notice that $r \ge 1$ implies that $x \in \G^r_1(x,M)$.
Furthermore, $\lim_{r \rightarrow \infty} \G^r_1(x,M) = \TY$ for all
$x \in \TY \setminus \Y_3$ and so we define $\G^{r=\infty}_1(x,M) =
\TY$ for all such $x$. For $x \in \Y_3$, $\G^r_1(x,M)= \{x\}$ for
all $r \in [1,\infty]$.

Notice that $\G^r_1(x,M_C)$ is a convex hexagon for all $r \ge 2$
and $x \in \TY\setminus \Y_3$, (since for such an $x$,
$\G^r_1(x,M_C)$ is bounded by $\xi_j(r,x)$ and $e_j$ for all $j \in
\{1,2,3\}$, see also Figure \ref{fig:ProxMapDef2},) else it is
either a \emph{convex hexagon} or a \emph{non-convex} but \emph{star-shaped
polygon} depending on the location of $x$ and the value of $r$.

\subsection{The $\G_1$-Region of a Multiple Point Data Set for $\NPE^r$}
So far, we have described the $\G_1$-region for a point in $x \in
\TY$. For a set $\X_n$ of size $n$ in $\TY$, the region
$\G^r_1(\X_n,M)$ can be specified by the edge extrema only.
The (closest) \emph{edge extrema} of a set $B$ in $\TY$ are the
points closest to the edges of $\TY$, denoted $x_{e_j}$ for $j\in
\{1,2,3\}$; that is, $x_{e_j} \in \arginf_{x \in B}d(x,e_j)$.
Note that if $B=\X_n$ is a set of iid random variables of size $n$
from $F$ then the edge extrema, denoted $X_{e_j}(n)$, are random
variables. Below, we show that the edge extrema are the active
points in defining $\G^r_1(\X_n,M)$.

\begin{proposition}
\label{prop:SMA-edge-extrema}
Let $B$ be any set of $n$ distinct points in $\TY$.
For $r$-factor proportional-edge proximity maps with $M$-vertex regions,
$\G^r_1\left(B,M \right)=\cap_{k=1}^3\,\G^r_1\left(x_{e_k},M \right)$.
\end{proposition}
{\bfseries Proof:}
Given $B = \{x_1,\ldots,x_n\}$ in $\TY$.
Note that
$$\G^r_1(B,M) \cap R_M(\y_j)=\bigl[\cap_{i=1}^n\,\G^r_1(x_i,M)\bigr]\cap R_M(\y_j),$$
but by definition  $x_{e_j} \in \argmax_{x\in B}d(\y_j,\xi_j(r,x))$,
so
 \begin{equation}
 \label{eqn:G1-in3}
 \G^r_1(B,M) \cap
R_M(\y_j)=\G^r_1(x_{e_j},M)\cap R_M(\y_j) \text{ for }j\in
\{1,2,3\}.
 \end{equation}
Furthermore, $\G^r_1(B,M) = \cup_{j=1}^3 \bigl[ \G^r_1(x_{e_j},M)
\cap R_M(\y_j) \bigr]$, and
 \begin{equation}
 \label{eqn:G1-each-in3}
 \G^r_1(x_{e_j},M) \cap
R_M(\y_j)=\cap_{k=1}^3 \bigl[ \G^r_1(x_{e_k},M) \cap R_M(\y_j)
\bigr] \text{ for }j\in \{1,2,3\}.
 \end{equation}
Combining these two results in Equations (\ref{eqn:G1-in3}) and (\ref{eqn:G1-each-in3}),
we obtain $\G^r_1(B,M)=\cap_{k=1}^3\,\G^r_1(x_{e_k},M)$. $\blacksquare$

\begin{figure}
\begin{center}
\psfrag{A}{}
 \psfrag{B}{}
\psfrag{C}{\scriptsize{$\y_3$}}
\psfrag{IC}{\scriptsize{$M_{I}$}}
 \psfrag{x}{}
 \psfrag{0}{\scriptsize{$\y_1$}}
 \psfrag{1}{\scriptsize{$\y_2$}}
\psfrag{x1}{\scriptsize{$x_{e_1}$} }
\psfrag{x2}{\scriptsize{$x_{e_2}$} }
\psfrag{x3}{\scriptsize{$x_{e_3}$} }
 \epsfig{figure=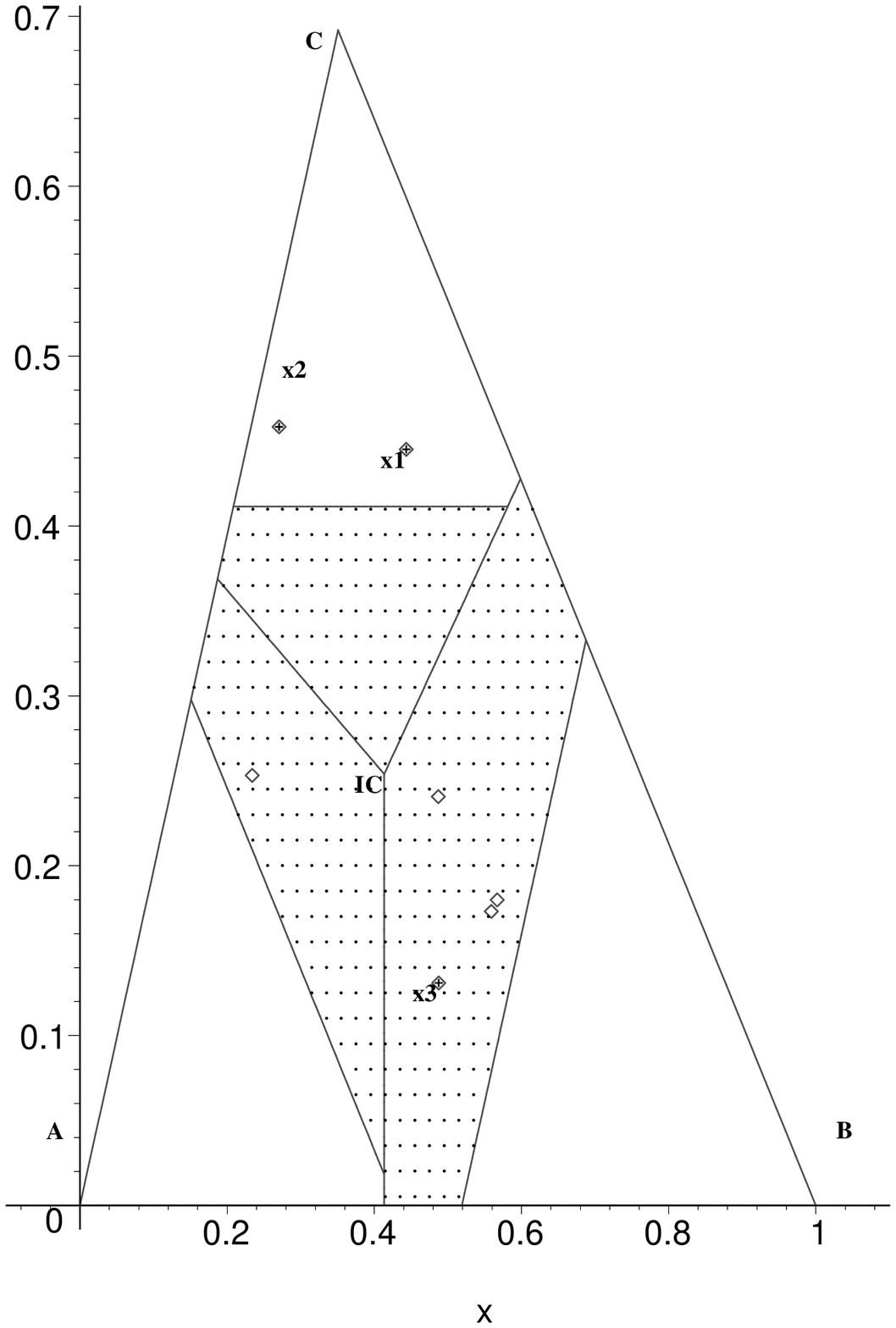, height=140pt,width=200pt}
\psfrag{(x1,y1)}{\scriptsize{$x_{e_1}$} }
\psfrag{(x2,y2)}{\scriptsize{$x_{e_2}$} }
\psfrag{(x3,y3)}{\scriptsize{$x_{e_3}$} }
\psfrag{CC}{\scriptsize{$M_{CC}$}}
 \epsfig{figure=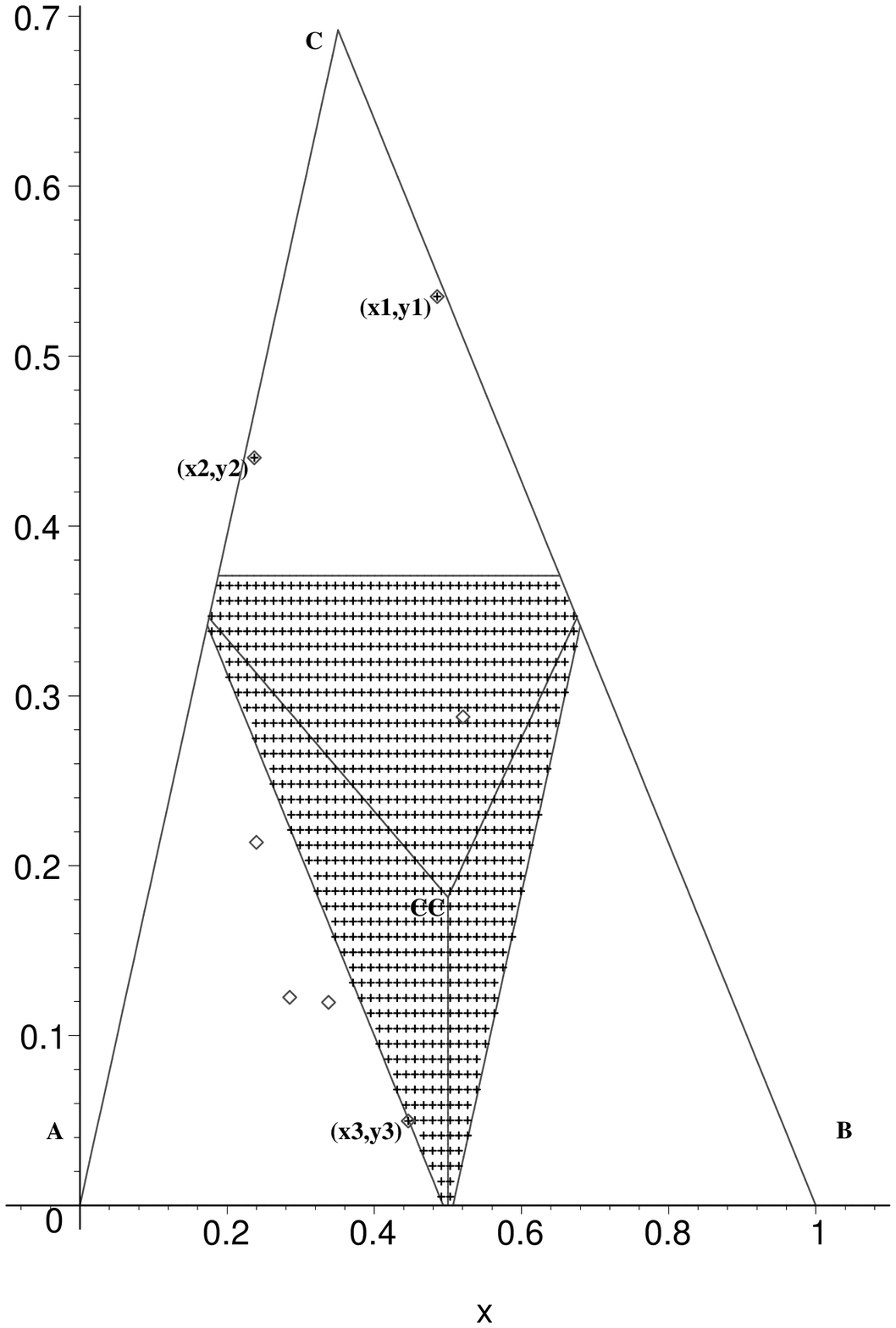, height=140pt,width=200pt}
\end{center}
\caption{
 \label{fig:G1-regions-n>1}
The $\G_1$-regions (the hatched regions) for $r=2$ with seven $\X$ points
iid $\U(\TY)$ where vertex regions constructed with incenter $M_I$
(left) and circumcenter $M_{CC}$ (right) with orthogonal projection.
}
\end{figure}

From the above proposition, we see that the $\G_1$-region for $B$ as
in proposition can also be written as the union of three regions of
the form
$$\G^r_1(B,M)\cap R_M(\y_j)=\{z \in R_M(\y_j):\;
d(\y_j, \ell(\y_j,z)) \ge d(\y_j,\xi_j(r,x_{e_j}))\} \text{ for } j\in \{1,2,3\}.$$

See Figure \ref{fig:G1-regions-n>1} for $\G_1$-region for $r=2$ with
seven $\X$ points iid $\U(\TY)$.
In the left figure, vertex regions are
based on incenter, while in the right figure, on circumcenter with
orthogonal projections to the edges. In either case $\X_n \cap
\G^{r=2}_1(\X_n,M)$ is nonempty, hence $\g_n(2,M)=1$.

Below, we demonstrate that
edge extrema are distinct with probability 1 as $n \rightarrow
\infty$. Hence in the limit three distinct points suffice to determine
the $\G_1$-region.
\begin{theorem}
\label{thm:distinct-edge-ext}
Let $\X_n$ be a set of iid random variables from $\UT$ and
let $E_{c,3}(n)$ be the event that (closest) edge extrema are distinct.
Then $P(E_{c,3}(n)) \rightarrow 1$ as $n \rightarrow \infty$.
\end{theorem}

We can also define the regions associated with $\g(\X_n,N)=k$ for $k
\le n$ called \emph{$\G_k$-region} for proximity map
$N_{\Y_3}(\cdot)$ and set $B \subseteq \Omega$ for $k=1,\ldots,n$
(see \cite{ceyhan:2004d}).

\section{The Asymptotic Distribution of $\g_n(r,M)$}
\label{sec:asy-gam-NYr} In this section, we first present a finite
sample upper bound for $\g_n(r,M)$, then present the degenerate
cases, and the nondegenerate case of the asymptotic distribution of
$\g_n(r,M)$ given in Equation (\ref{eqn:asymptotic-NYr}).

\subsection{An Upper Bound for $\g_n(r,M)$}
 Recall that by definition, $\g(\X_n,N) \le
n$. We will seek an a.s. least upper bound for $\g(\X_n,N)$. Let
$\X_n$ be a set of iid random variables from $F$ on $\TY$ and let
$\g(\X_n,N)$ be the domination number for the PCD based on a
proximity map $N$. Denote the general a.s. least upper bound for
$\g(\X_n,N)$ that works for all $n \ge 1$ and is independent of $n$
(which is called \emph{$\kappa$-value} in \cite{ceyhan:2004d}) as
$\kappa(N) := \min \{k: \g(\X_n, N) \le k \text{ a.s. for all } n
\ge 1\}$.

In $\R$ with $\Y_2=\{0,1\}$, for $\X_n$ a set of iid random variables from $\U(0,1)$,
$\g(\X_n,\NS) \le 2$ with equality holding with positive probability.
Hence $\kappa(\NS)=2$.

\begin{theorem}
\label{thm:kappaNYr=3}
Let $\X_n$ be a set of iid random variables from $\UT$ and $M \in \R^2 \setminus \Y_3$.
Then $\kappa\left(\NPE^r \right)=3$ for $\NPE^r(\cdot,M)$.
\end{theorem}
{\bfseries Proof:}
For $\NPE^r(\cdot,M)$, pick the point closest to edge $e_j$
in vertex region $R_M(\y_j)$;
that is, pick
$U_j \in \argmin_{X \in \X_n \cap R_M(\y_j)} d(X,e_j)
=\argmax_{X \in \X_n \cap R_M(\y_j)} d(\ell(\y,X),\y_j)$
in the vertex region for which $X_n \cap R_M(\y_j) \not=\emptyset$
for $j \in \{1,2,3\}$ (note that as $n \rightarrow \infty$,
$U_j$ is unique a.s. for each $j$,
since $X$ is from $\UT$).
Then $\X_n \cap R_M(\y_j) \subset \NPE^r(U_j,M)$.
Hence $\X_n \subset \cup_{j=1}^3 \NPE^r(U_j,M)$.
 So $\g_n(r,M_C) \le 3$ with equality holding with positive probability.
Thus $\kappa\left(\NPE^r \right)=3$.
$\blacksquare$

Below is a general result for the limiting distribution of
$\g(\X_n,N)$ for $\X_n$ from a very broad family of distributions
and for general $N(\cdot)$.

\begin{lemma}
\label{lem:RS-nonempty} Let $\RS(N)$ be the superset region for the
proximity map $N(\cdot)$ and $\X_n$ be a set of iid random variables
from $F$ with $P_F(X \in \RS(N))>0$. Then $\lim_{n \rightarrow
\infty}P_F(\g(\X_n,N)=1)=1$.
\end{lemma}
{\bfseries Proof:}
Suppose $P_F(X \in \RS(N))>0$.
Recall that for any $x \in \RS(N)$, we have $N(x) = \Omega$,
so $\X_n \subseteq N(x)$, hence if $\X_n \cap \RS(N) \not= \emptyset$
then $\g(\X_n,N)=1$.
Then $P(\X_n \cap \RS(N)\not= \emptyset) \le P(\g(\X_n,N)=1)$.
But
$P\bigl(\X_n \cap \RS(N)\not= \emptyset \bigr)=
1-P\bigl( \X_n \cap \RS(N)=\emptyset \bigr)=
1-\bigl[1-P_F\bigl(X \in \RS(N)\bigr)\bigr]^n
\rightarrow 1 \text{ as }n \rightarrow \infty,$
since $P_F\bigl(X \in \RS(N)\bigr)>0$.
Hence $\lim_{n \rightarrow \infty}P(\g(\X_n,N)=1) =1$.
$\blacksquare$

\begin{remark}
\label{rem:rate-of-convergence}
In particular, for $F=\UT$,
the inequality $P_F(X \in \RS(N))>0$ holds iff $A(\RS(N))>0$,
then $P(\X_n \cap \RS(N)\not= \emptyset)\rightarrow 1$. $\square$
\end{remark}
For $\Y_2=\{0,1\} \subset \R$, $\RS(\NS)=\{1/2\}$, so Lemma
\ref{lem:RS-nonempty} does not apply to $N_S$ in $\R$.

Recall that $\kappa\left( \NPE^r \right) = 3$, then
$$1 \le \E\left[ \g_n(r,M) \right] \le 3 \text{    and    } 0 \le \Var\left[ \g_n(r,M) \right] \le 9/4.$$

Furthermore, there is a stochastic ordering for $\g_n(r,M)$.

\begin{theorem}
\label{thm:stochastic-r1<r2} Suppose $\X_n$ is a set of iid random
variables from a continuous distribution $F$ on $\TY$. Then for
$r_1<r_2$, we have $\g_n(r_2,M) \le^{ST} \g\left(\X_n, \NPE^{r_1},M
\right)$.
\end{theorem}

{\bfseries Proof:}
Suppose $r_1<r_2$.
Then
$P\left(\g_n(r_2,M) \le 1\right)>P\left(\g_n(r_1,M)\le 1\right)$
since $\G^{r_1}_1(\X_n,M) \subsetneq \G^{r_2}_1(\X_n,M)$
for any realization of $\X_n$ and by a similar argument
$P\bigl(\g_n(r_2,M) \le 2\bigr)>P\left(\g_n(r_1,M) \le 2\right)$
so
$P\left(\g_n(r_2,M) \le 3\right) = P\left(\g_n(r_1,M) \le 3\right).$
Hence the desired result follows. $\blacksquare$

\subsection{Geometry Invariance}
We present a ``geometry invariance" result for $\NPE^r(\cdot,M)$
where $M$-vertex regions are constructed using the lines joining
$\Y_3$ to $M$, rather than the orthogonal projections from $M$ to
the edges. This invariance property will simplify the notation in
our subsequent analysis by allowing us to consider the special case
of the equilateral triangle.

\begin{theorem}
\label{thm:geo-inv-NYr} (Geometry Invariance Property)
Suppose $\X_n$ is a set of iid random variables from $\UT$.
Then for any $r \in [1,\infty]$ the distribution of $\g_n(r,M)$ is
independent of $\Y_3$ and hence the geometry of $\TY$.
\end{theorem}

{\bfseries Proof:}
Suppose $X \sim \U(T(\Y))$.
A composition of translation, rotation,
reflections, and scaling will take any given triangle
$T(\Y)=T(\y_1,\y_2,\y_3)$ to the basic triangle
$T_b=T((0,0),(1,0),(c_1,c_2))$ with $0 < c_1 \le 1/2$, $c_2>0$, and
$(1-c_1)^2+c_2^2 \le 1$.
Furthermore, when $X$ is also transformed in the same manner, say to $X'$,
then $X'$ is uniform on $T_b$, i.e., $X' \sim \mathcal U(T_b)$.
The transformation
$\phi_e:\R^2 \rightarrow \R^2$ given by $\phi_e(u,v)=\left(
u+\frac{1-2\,c_1}{\sqrt{3}}\,v, \frac{\sqrt{3}}{2\,c_2}\,v \right)$
takes $T_b$ to the equilateral triangle
$T_e=\bigl((0,0),(1,0),(1/2,\sqrt{3}/2)\bigr)$.
Investigation of the
Jacobian shows that $\phi_e$ also preserves uniformity.
That
is, $\phi_e(X') \sim \mathcal U(T_e)$.
Furthermore,
the composition of $\phi_e$, with the scaling and rigid body
transformations, maps
     the boundary of the original triangle, $T_o$,
  to the boundary of the equilateral triangle, $T_e$,
     the lines joining $M$ to $\y_j$ in $T_b$
  to the lines joining $\phi_e(M)$ to $\phi_e(\y_j)$ in $T_e$,
and  lines parallel to the edges of $T_o$
  to lines parallel to the edges of $T_e$.
Since the distribution of $\g_n(r,M)$ involves only probability
content of unions and intersections of regions bounded by precisely
such lines and the probability content of such regions is preserved
since uniformity is preserved; the desired result follows.
$\blacksquare$

Note that geometry invariance of $\g\left(\X_n, \NPE^{r=\infty},M
\right)$ also follows trivially, since for $r=\infty$,
we have $\g_n(r=\infty,M)=1$ a.s. for all $\X_n$ from any $F$ with support
in $\TY \setminus \Y_3$.

Based on Theorem \ref{thm:geo-inv-NYr} we may assume that $\TY$ is a
standard equilateral triangle with \\$\Y_3= \left\{ (0,0),(1,0),
\left( 1/2,\sqrt{3}/2 \right) \right\}$ for $\NPE^r(\cdot,M)$ with
$M$-vertex regions.

Notice that, we proved the geometry invariance property for $\NPE^r$
where $M$-vertex regions are defined with the lines joining $\Y_3$ to $M$.
On the other hand, if we use the orthogonal projections from
$M$ to the edges, the vertex regions, hence $\NPE^r$ will depend on
the geometry of the triangle.
That is, the orthogonal projections
from $M$ to the edges will not be mapped to the orthogonal
projections in the standard equilateral triangle.
Hence with the
choice of the former type of $M$-vertex regions, it suffices to work
on the standard equilateral triangle.
On the other hand, with the orthogonal projections, the exact and
asymptotic distribution of $\g_n$ will depend on $c_1,c_2$, so one
needs to do the calculations for each possible combination of
$c_1,c_2$.

\subsection{The Degenerate Case with $\g_n(r,M) \stackrel{p}{\rightarrow} 1$}

Below, we prove that $\g_n(r,M)$ is degenerate in the limit for
$r>3/2$.
\begin{theorem}
\label{thm:M-notin-Tr} Suppose $\X_n$ is a set of iid random
variables from a continuous distribution $F$ on $\TY$. If $M \not\in
\Tr$ (see Figure \ref{fig:Tr-RS-NDA-CC} and Equation
(\ref{eqn:T^r-def}) for $\Tr$), then $\lim_{n \rightarrow
\infty}P\left( \g_n(r,M)=1 \right) =1$ for all $M \in \R^2 \setminus
\Y_3$.
\end{theorem}
{\bfseries Proof:} Suppose $M \notin \Tr$. Then $\RS\left( \NPE^r,M
\right)$ is nonempty with positive area. Hence the result follows by
Lemma \ref{lem:RS-nonempty}. $\blacksquare$

\begin{corollary}
\label{cor:M-notin-Tr} Suppose $\X_n$ is a set of iid random
variables from a continuous distribution $F$ on $\TY$. Then for
$r>3/2$,
 $\lim_{n \rightarrow \infty}P\left( \g_n(r,M)=1 \right) =1$
for all $M \in R^2\setminus \Y_3$.
\end{corollary}
{\bfseries Proof:}
For $r>3/2$,  $\Tr=\emptyset$, so $M \not\in \Tr$.
Hence the result follows by Theorem \ref{thm:M-notin-Tr}.
$\blacksquare$

We estimate the distribution of $\g_n(r,M)$ with $r=2$ and $M=M_C$
for various $n$ empirically. In Table
\ref{tab:numerical-gamma-NYr=2,5/4} (left), we present the empirical
estimates of $\g_n(r,M)$ with $n=10,\,20,\,30,\,50,\,100$ based on
$1000$ Monte Carlo replicates in $T_e$. Observe that the empirical
estimates are in agreement with the asymptotic distribution given in
Corollary \ref{cor:M-notin-Tr}.

\begin{table}[ht]
\begin{center}
\begin{tabular}{cc}

\begin{tabular}{|c|c|c|c|c|c|}
\hline
$ k \diagdown n$  & 10 & 20 & 30 & 50 & 100\\
\hline
1 & 961 & 1000 & 1000 & 1000 & 1000 \\
\hline
2 & 34 & 0 & 0 & 0 & 0  \\
\hline
3 & 5 & 0 & 0 & 0 & 0 \\
\hline
\end{tabular}

&

\begin{tabular}{|c|c|c|c|c|c|}
\hline
$ k \diagdown n$  & 10 & 20 & 30 & 50 & 100\\
\hline
1 & 9 & 0 & 0 & 0 & 0 \\
\hline
2 & 293 & 110 & 30 & 8 & 0  \\
\hline
3 & 698 & 890 & 970 & 992 & 1000 \\
\hline
\end{tabular}

\end{tabular}
\end{center}
\caption{ \label{tab:numerical-gamma-NYr=2,5/4} The number of
$\g_n(r,M)=k$ out of $N=1000$ Monte Carlo replicates with $M=M_C$
and $r=2$ (left) and $r=5/4$ (right).}
\end{table}


The asymptotic distribution of $\g_n(r,M)$ for $r <3/2$ depends on
the relative position of $M$ with respect to the triangle $\Tr$.

\subsection{The Degenerate Case with $\g_n(r,M) \stackrel{p}{\rightarrow} 3$}

\begin{theorem}
\label{thm:gnNYr=3}
Suppose $\X_n$ is a set of iid random variables from a continuous distribution $F$
on $\TY$.
If $M \in (\Tr)^o$, then $P\left( \g_n(r,M)=3 \right) \rightarrow 1$
as $n \rightarrow \infty$.
\end{theorem}

We estimate the distribution of $\g_n(r,M)$ with $r=5/4$ and $M=M_C$
for various $n$ values empirically. In Table
\ref{tab:numerical-gamma-NYr=2,5/4} (right), we present the
empirical estimates of $\g_n(r,M)$ with $n=10,\,20,\,30,\,50,\,100$
based on $1000$ Monte Carlo replicates in $T_e$. Observe that the
empirical estimates are in agreement with our result in Theorem
\ref{thm:gnNYr=3}.

\begin{theorem}
\label{thm:MinboundaryTr-g>1} Suppose $\X_n$ is a set of iid random
variables from $\UT$. If $M \in \partial(\Tr)$, then $P\left(
\g_n(r,M)>1 \right) \rightarrow 1$ as $n \rightarrow \infty$.
\end{theorem}

For $M \in \partial(\Tr)$, there are two separate cases:
\begin{itemize}
\item[(i)]
$M \in \partial(\Tr) \setminus \{t_1(r),\,t_2(r),\,t_3(r)\}$ where $t_j(r)$ with
$j \in \{1,2,3\}$ are the vertices of $\Tr$ whose explicit forms are
given in Equation (\ref{eqn:T^r-def}).
\item[(ii)]
$M \in \{t_1(r),\,t_2(r),\,t_3(r)\}$.
\end{itemize}

\begin{theorem}
\label{thm:g_nNYr=3-for-r<3/2}
Suppose $\X_n$ is a set of iid random variables from  $\UT$.
If $M \in \partial(\Tr)\setminus \{t_1(r),\,t_2(r),\,t_3(r)\}$, then
$P\left( \g_n(r,M)=3 \right) \rightarrow 1$ as $n \rightarrow
\infty$.
\end{theorem}

We estimate the distribution of $\g_n(r,M)$ with $r=5/4$ and
$M=\left( 3/5,\sqrt{3}/10\right)\in \partial(\Tr)\setminus
\{t_1(r),\,t_2(r),\,t_3(r)\}$ for various $n$ empirically. In Table
\ref{tab:numerical-gamma-NYr=5/4-1} we present empirical estimates
of $\g_n(r,M)$ with $n=10,\,20,\,30,\,50,\,100,\,500,\\
\,1000,\,2000$
based on $1000$ Monte Carlo replicates in $T_e$. Observe that the
empirical estimates are in agreement with our result in Theorem
\ref{thm:g_nNYr=3-for-r<3/2}.

\begin{table}[ht]
\begin{center}
\begin{tabular}{|c|c|c|c|c|c|c|c|c|}
\hline
$ k \diagdown n$  & 10 & 20 & 30 & 50 & 100 & 500 & 1000 & 2000\\
\hline
1 & 118 & 60 & 51 & 39 & 15 & 1 & 2 & 1 \\
\hline
2 & 462 & 409 & 361 & 299 & 258 & 100 & 57 & 29  \\
\hline
3 & 420 & 531 & 588 & 662 & 727 & 899 & 941 & 970\\
\hline
\end{tabular}
\end{center}
\caption{\label{tab:numerical-gamma-NYr=5/4-1} The number of
$\g_n(r,M)=k$ out of $N=1000$ Monte Carlo replicates with $r=5/4$
and $M=\left( 3/5,\sqrt{3}/10\right)$.}
\end{table}

\subsection{The Nondegenerate Case}

\begin{theorem}
\label{thm:g_nNYr=2-for-r-vertex-of-Tr}
Suppose $\X_n$ is a set of iid random variables from $\UT$.
If $M \in \{t_1(r),\,t_2(r),\,t_3(r)\}$, then $P\left( \g_n(r,M)=2
\right) \rightarrow p_r$ as $n \rightarrow \infty$ where $p_r \in (0,1)$ is
provided in Equation (\ref{eqn:p_r-form}) but only
numerically computable.
\end{theorem}

For example, $p_{r=5/4} \approx .6514$ and $p_{r=\sqrt{2}} \approx .4826$.

So the asymptotic distribution of $\g_n(r,M)$ with $r \in [1,3/2)$
and $M \in \{t_1(r),t_2(r),t_3(r)\}$ is given by
\begin{equation}
\label{eqn:Asydist-NYr} \g_n(r,M) \sim 2+\BER(1-p_r).
\end{equation}

We estimate the distribution of $\g_n(r,M)$ with $r=5/4$ and
$M=\left( 7/10,\sqrt{3}/10 \right)$ for various $n$ empirically. In
Table \ref{tab:numerical-gamma-NYr=5/4-2}, we present the empirical
estimates of $\g_n(r,M)$ with
$n=10,\,20,\,30,\,50,\,100,\,500,\,1000,\,2000$ based on $1000$
Monte Carlo replicates in $T_e$. Observe that the empirical
estimates are in agreement with our result $p_{r=5/4} \approx
.6514$.

\begin{table}[ht]
\begin{center}
\begin{tabular}{|c|c|c|c|c|c|c|c|c|}
\hline
$ k \diagdown n$  & 10 & 20 & 30 & 50 & 100 & 500 & 1000 & 2000\\
\hline
1 & 174 & 118 & 82 & 61 & 22 & 5 & 1 & 1\\
\hline
2 & 532 & 526 & 548 & 561 & 611 & 617 & 633 & 649  \\
\hline
3 & 294 & 356 & 370 & 378 & 367 & 378 & 366 & 350\\
\hline
\end{tabular}
\end{center}
\caption{\label{tab:numerical-gamma-NYr=5/4-2} The number of
$\g_n(r,M)=k$ out of $N=1000$ Monte Carlo replicates with $r=5/4$
and $M=\left(7/10,\sqrt{3}/10 \right)$.}
\end{table}

\begin{remark}
For $r=3/2$, as $n \rightarrow \infty$, $P\left( \g_n(r,M_C)
> 1 \right) \rightarrow 1$ at rate $O\left( n^{-1} \right)$.
$\square$
\end{remark}

\begin{theorem}
\label{thm:asydist-gamma-NY3/2}
Suppose $\X_n$ is a set of iid random variables from $\UT$.
Then for $r=3/2$, as $n \rightarrow \infty$,
\begin{equation}
\label{eqn:Asydist-NY3/2} \g_n(3/2,M_C) \sim 2+\BER(p \approx .2487)
\end{equation}
\end{theorem}

For the proof of Theorem \ref{thm:asydist-gamma-NY3/2},
see \cite{ceyhan:2004e, ceyhan:2005e}.

Using Theorem \ref{thm:asydist-gamma-NY3/2},
\begin{equation}
\label{eqn:Asymom-NY3/2-1} \lim_{n \rightarrow \infty}\E\left[
\g_n(3/2,M_C) \right]= 3-p_{3/2} \approx 2.2587
\end{equation}
and
\begin{equation}
\label{eqn:Asymom-NY3/2-2} \lim_{n \rightarrow \infty}\Var\left[
\g_n(3/2,M_C) \right]= 6+p_{3/2}-p_{3/2}^2 \approx .1917.
\end{equation}
Indeed, the finite sample distribution of $\g_n(3/2,M_C)$ hence the
finite sample mean and variance can also be obtained by numerical
methods.

We also estimate the distribution of $\g_n(3/2,M_C)$ for various $n$
values empirically.
The empirical estimates for
$n=10,\,20,30,\,50,\,100,\,500,\,1000,\,2000$ based on $1000$ Monte
Carlo replicates are given in Table \ref{tab:numerical-gamma-NY3/2}.
estimates are in agreement with our result $p_{r=3/2} \approx
.7413$.

\begin{table}[ht]
\begin{center}
\begin{tabular}{|c|c|c|c|c|c|c|c|c|}
\hline
$ k \diagdown n$  & 10 & 20 & 30 &  50 & 100 & 500 & 1000 & 2000\\
\hline
1  & 151 & 82 & 61 & 50 & 27 & 2 & 3 & 1\\
\hline
2 &  602 & 636 & 688 & 693 & 718 & 753  & 729 & 749\\
\hline
3  & 247 & 282 & 251 & 257 & 255 & 245 & 268 & 250 \\
\hline
\end{tabular}
\end{center}
\caption{\label{tab:numerical-gamma-NY3/2} The number of
$\g_n(3/2,M_C)=k$ out of $N=1000$ Monte Carlo replicates.}
\end{table}

\subsection{Distribution of the $\g_n(r,M)$ in Multiple Triangles}
\begin{figure}[ht]
\centering
 \rotatebox{-90}{ \resizebox{3. in}{!}{\includegraphics{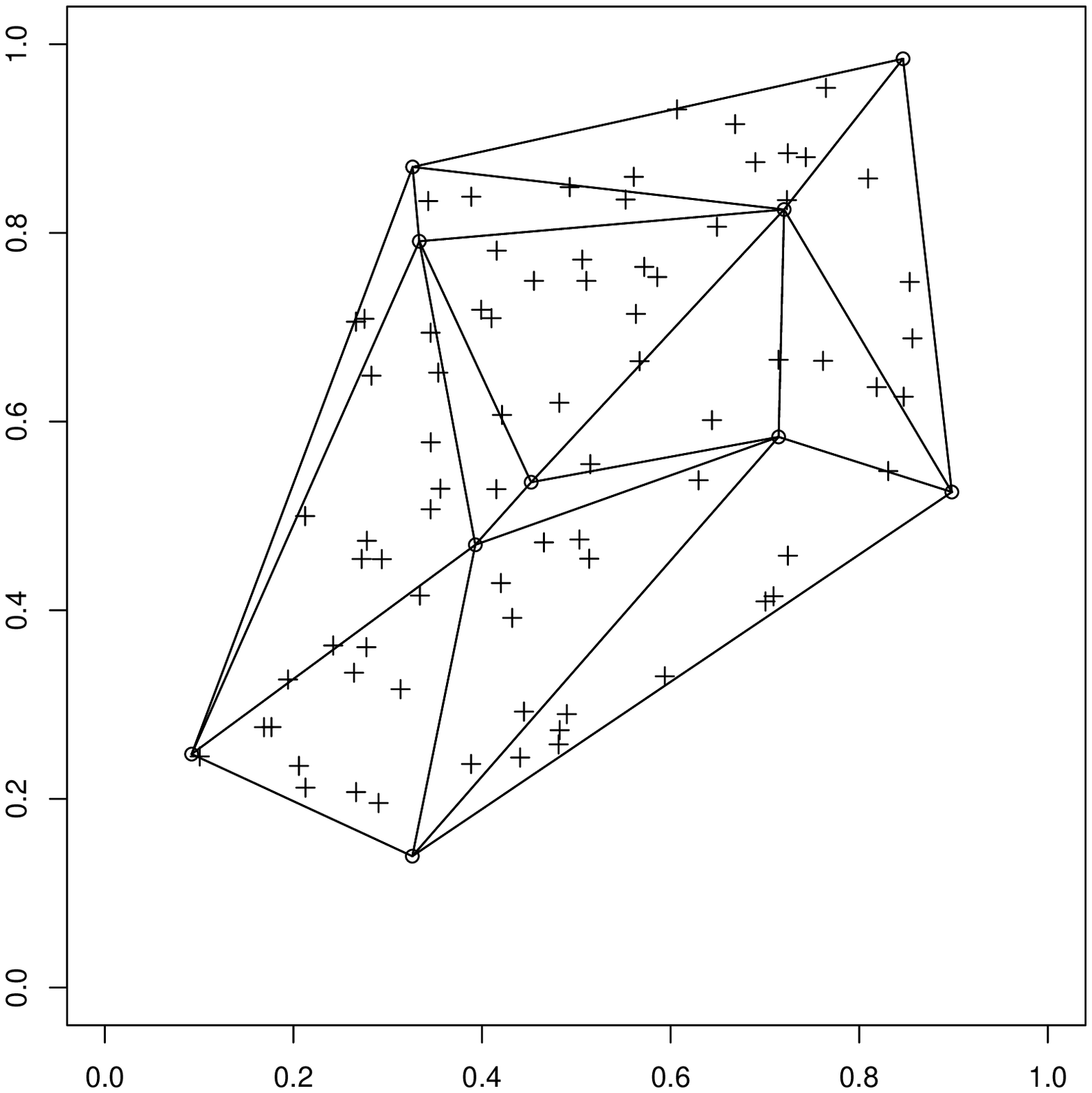}}}
 \rotatebox{-90}{ \resizebox{3. in}{!}{\includegraphics{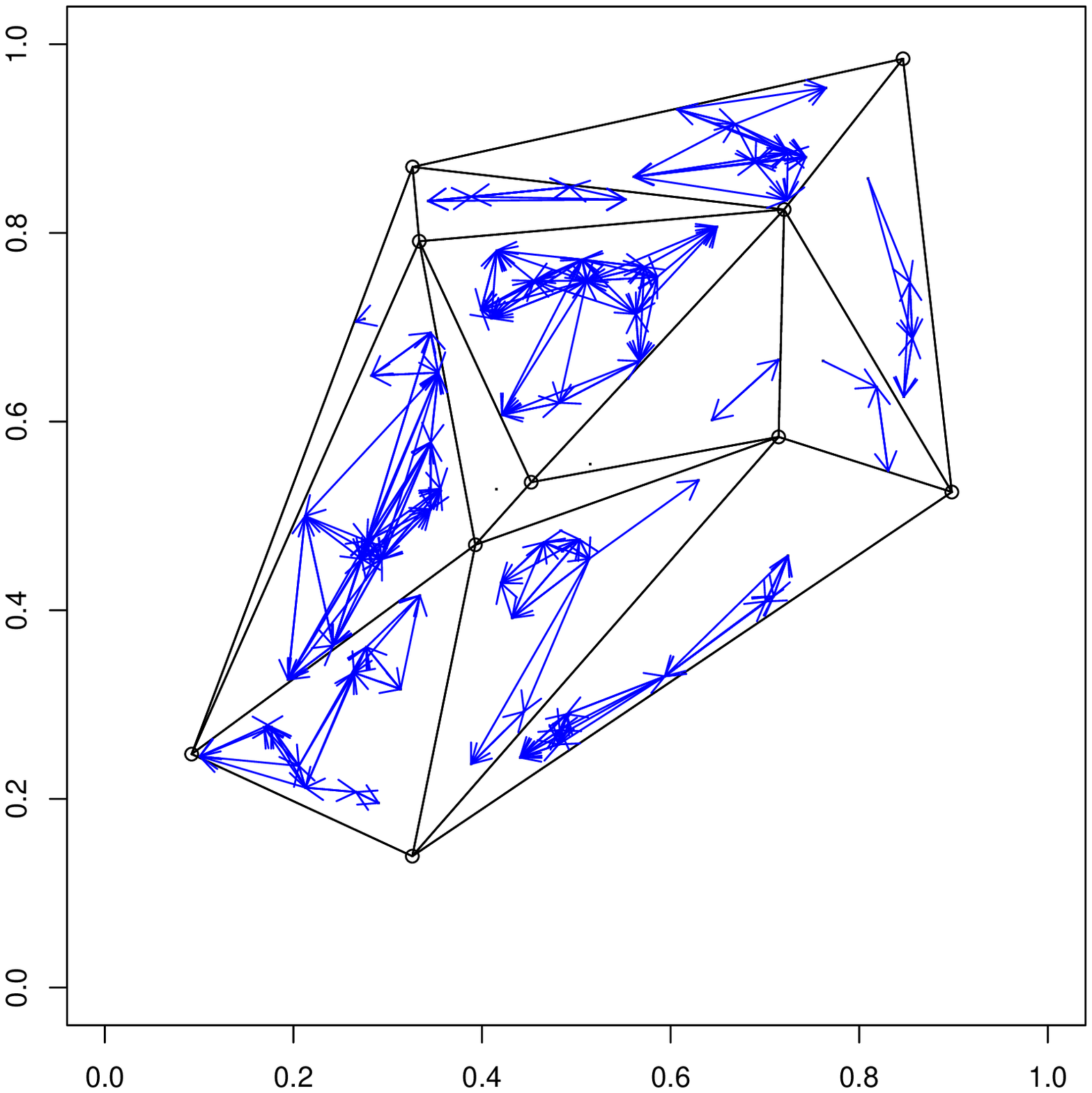}}}
 \caption{
 \label{fig:multi-tri-arcs}
 The 77 $\X$ points (crosses) in the convex hull of $\Y$ points (circles)
 given in Figure \ref{fig:deltri} (left) and the corresponding arcs (right) of
 $r$-factor proportional edge PCD with $r=3/2$ and $M=M_C$. }
\end{figure}

So far we have worked with data in one Delaunay triangle, i.e.,
$m=3$ or $J_3=1$. In this section, we present the asymptotic
distribution of the domination number of $r$-factor PCDs in multiple
Delaunay triangles. Suppose $\Y_m=\{\y_1,\y_2,\ldots,\y_m\} \subset
\R^2$ be a set of $m$ points in general position with $m>3$ and no
more than 3 points are cocircular. Then there are $J_m>1$ Delaunay
triangles each of which is denoted as $\T_j$. Let $M^j$ be the point
in $\T_j$ that corresponds to $M$ in $T_e$, $\Tr^j$ be the triangle
that corresponds to $\Tr$ in $T_e$, and $t_i^j(r)$ be the vertices
of $\Tr^j$ that correspond to $t_i(r)$ in $T_e$ for $i \in
\{1,2,3\}$. Moreover, let $n_j:=|\X_n \cap \T_j|$, the number of $X$
points in Delaunay triangle $\T_j$. For $\X_n \subset \C_H(\Y_m)$,
let $\g_{n_j}(r,M^j)$ be the domination number of the digraph
induced by vertices of $\T_j$ and $\X_n \cap \T_j$. Then the
domination number of the $r$-factor PCD in $J_m$ triangles is
$$\g_n(r,M,J_m)=\sum_{j=1}^{J_m}\g_{n_j}(r,M^j).$$
See Figure \ref{fig:multi-tri-arcs} (left) for the 77 $\X$ points
that are in $\C_H(\Y_m)$ out of the 200 $\X$ points plotted in
Figure \ref{fig:deltri}. Observe that 10 $\Y$ points yield
$J_{10}=13$ Delaunay triangles. In Figure \ref{fig:multi-tri-arcs}
(right) are the corresponding arcs for $M=M_C$ and $r=3/2$. The
corresponding $\g_n=22$. Suppose $\X_n$ is a set of iid random
variables from $\U(\C_H(\Y_m))$, the uniform distribution on convex
hull of $\Y_m$ and we construct the $r$-factor PCDs using the points
$M^j$ that correspond to $M$ in $T_e$. Then for fixed $m$ (or fixed
$J_m$), as $n \rightarrow \infty$, so does each $n_j$. Furthermore,
as $n \rightarrow \infty$, each component $\g_{n_j}(r,M^j)$ become
independent. Therefore using Equation (\ref{eqn:asymptotic-NYr}), we
can obtain the asymptotic distribution of $\g_n(r,M,J_m)$. As $n
\rightarrow \infty$, for fixed $J_m$,\\

\noindent
 \begin{equation}
 \label{eqn:asymptotic-NYr-Jm}
 \g_n(r,M,J_m) \sim
\left\lbrace
 \begin{array}{ll}
       2\,J_m+\BIN(J_m,1-p_r),           & \text{for $M^j \in \{t^j_1(r),t^j_2(r),t^j_3(r)\}$ and $r \in [1,3/2]$,}\\
       J_m,           & \text{for $r>3/2$,}\\
       3\,J_m,           & \text{for $M \in \Tr^j\setminus \{t^j_1(r),t^j_2(r),t^j_3(r)\}$ and $r \in [1,3/2)$,}\\

\end{array} \right.
\end{equation}
\noindent where $\BIN(n,p)$ stands for binomial distribution with
$n$ trials and probability of success $p$, for $r \in [1,3/2)$ and
$M \in \{t_1(r),t_2(r),t_3(r)\}$, $p_r$ is given in Equation
\ref{eqn:asymptotic-NYr} and for $r=3/2$ and $M=M_C$, $p_r \approx
.7413$ (see Equation (\ref{eqn:Asydist-NY3/2})).

\subsection{Extension of $\NPE^r$ to Higher Dimensions}
\label{sec:NYr-higher-D}
The extension to $\R^d$ for $d > 2$ with $M=M_C$ is provided in \cite{ceyhan:2005e},
but the extension for general $M$ is similar.

Let $\g_n(r,M,d):=\g(\X_n,\NPE^r,M,d)$ be the domination number of the PCD based on
the extension of $\NPE^r(\cdot,M)$ to $\R^d$. Then it is easy to
see that $\g_n(r,M,3)$ is nondegenerate as $n \rightarrow \infty$ for $r=4/3$.
In $\R^d$, it can be seen that
$\g_n(r,M,d)$ is nondegenerate in the limit only when $r=(d+1)/d$.
Furthermore, for large $d$, asymptotic distribution of
$\g_n(r,M,d)$ is nondegenerate at values of $r$ closer to $1$.
Moreover, it can be shown that
$\lim_{n \rightarrow \infty}P \bigl( 2 \le \g_n(r=(d+1)/d,M,d ) \le d+1\bigr)=1$
and we conjecture the following.

\begin{conjecture}
Suppose $\X_n$ is  set of iid random variables from the uniform
distribution on a simplex in $\R^d$. Then the domination number
$\g_n(r,M)$ in the simplex satisfies
$$\lim_{n \rightarrow \infty}P\left( d \le \g_n((d+1)/d,M,d) \le d+1 \right)=1.$$
\end{conjecture}

For instance, with  $d=3$ we estimate the empirical distribution of
$\g(\X_n,4/3)$ for various $n$. The empirical estimates for
$n=10,\,20,\,30,\,40,\,50,\,100,\,200,\,500,\,1000,\,2000$ based on $1000$
Monte Carlo replicates for each $n$ are given in
Table \ref{tab:numerical-gamma-NY4/3-in-R3}.

\begin{table}[ht]
\begin{center}
\begin{tabular}{|c|c|c|c|c|c|c|c|c|c|c|}
\hline
$ k \diagdown n$  & 10 & 20 & 30 & 40 & 50 & 100 & 200 & 500 & 1000 & 2000\\
\hline
1  & 52 &  18 &   5 &   5  &  4  &  0  &  0  &  0 &   0  &   0\\
\hline
2 &  385 & 308 & 263 & 221 & 219 & 155 & 88 &  41 &  31 &   19 \\
\hline
3  & 348 &  455 &  557 & 609 & 621 & 725 & 773 & 831 & 845 &  862\\
\hline
4  & 215 &  219 &  175 & 165 & 156 & 120 & 139 & 128 & 124 &  119\\
\hline
\end{tabular}
\end{center}
\caption{\label{tab:numerical-gamma-NY4/3-in-R3} The number of
$\g_n(4/3,M_C)=k$ out of $N=1000$ Monte Carlo replicates.}
\end{table}

\section{Discussion}
The $r$-factor proportional-edge proximity catch digraphs (PCDs),
when compared to class cover catch digraphs (CCCDs), have some advantages.
The asymptotic distribution of the domination number
$\g_n(r,M)$ of the $r$-factor PCDs, unlike that of CCCDs, is mathematically tractable
(computable by numerical integration).
A minimum dominating set can be found in
polynomial time for $r$-factor PCDs in $\R^d$ for all $d\ge 1$,
but finding a minimum
dominating set is an NP-hard problem for CCCDs (except for $\R$).
These nice properties of $r$-factor PCDs are due to the geometry invariance of
distribution of $\g_n(r,M)$ for uniform data in triangles.

On the other hand,
CCCDs are easily extendable to higher dimensions
and are defined for all $\X_n \subset \R^d$,
while $r$-factor PCDs are only defined for
$\X_n \subset \C_H(\Y_m)$.
Furthermore, the CCCDs based on balls use proximity regions
that are defined by the obvious metric, while the PCDs
in general do not suggest a metric.
In particular, our $r$-factor PCDs are based on some sort of dissimilarity measure,
but no metric underlying this measure exists.

The finite sample distribution of $\g_n(r,M)$, although
computationally tedious, can be found by numerical methods, while
that of CCCDs can only be empirically estimated by Monte Carlo simulations.
Moreover, we had to introduce many auxiliary tools to compute the distribution
of  $\g_n(r,M)$ in $\R^2$.
Same tools will work in higher dimensions,
perhaps with more complicated geometry.

The $r$-factor PCDs have applications in classification and testing
spatial patterns of segregation or association.
The former can be performed building discriminant regions for classification
in a manner analogous to the procedure proposed in \cite{priebe:2003b};
and the latter can be performed by using the asymptotic
distribution of $\g_n(r,M)$ similar to the procedure used in \cite{ceyhan:2005e}.

\section*{Acknowledgements}
This work was partially by the
Defense Advanced Research Projects Agency
as administered by the Air Force Office of Scientific Research
under contract DOD F49620-99-1-0213 and by
Office of Naval Research Grant N00014-95-1-0777.
We also thank anonymous referees, whose constructive comments and
suggestions greatly improved the presentation and flow of this article.


\section*{Appendix}

{\small
First, we begin with a remark that introduces some
terminology which we will use for asymptotics throughout this
appendix.

\begin{remark}
\label{rem:asymp-accurate} Suppose $\X_n$ is a set of iid random
variables from $F$ with support $\mS(F) \subseteq \Omega$. If over a
sequence $\Omega_n \subseteq \Omega,\; n=1,2,3,\ldots$, $X$
restricted to $\Omega_n$, $X|_{\Omega_n}$, has distribution $F_n$
with $F_n(x)=F(x)/P_F(X \in \Omega_n)$ and $P_F(X \in \Omega_n)
\rightarrow 1$ as $n \rightarrow \infty$, then we call $F_n$ the
\emph{asymptotically accurate distribution} of $X$ and $\Omega_n$
the \emph{asymptotically accurate support} of $F$ . If $F$ has
density $f$, then $f_n=f(x)/P_F(X \in \Omega_n)$ is called the
\emph{asymptotically accurate pdf} of $X$. In both cases, if we
are concerned with asymptotic results, for simplicity we will,
respectively, use $F$ and $f$ for asymptotically accurate
distribution and pdf. Conditioning will be implied by stating
that $X \in \Omega_n$
 with probability 1, as $n \rightarrow \infty$
or for sufficiently large $n$. $\square$
\end{remark}

\subsubsection*{Proof of Theorem \ref{thm:distinct-edge-ext}}
 Without loss of generality, assume $\TY=T_b=T((0,0),(1,0),(c_1,c_2))$ Note
that the probability of edge extrema all being equal to each other
is $P(X_{e_1}(n)=X_{e_2}(n)=X_{e_3}(n))=\I(n=1)$. Let $E_{c,2}(n)$
be the event that there are only two distinct (closest) edge
extrema. Then for $n>1$,
$$P(E_{c,2}(n))=P(X_{e_1}(n)=X_{e_2}(n))+P(X_{e_1}(n)=X_{e_3}(n))+P(X_{e_2}(n)=X_{e_3}(n))$$
since the intersection of the events $\{X_{e_i}(n)=X_{e_j}(n)\}$ and
$\{X_{e_i}(n)=X_{e_k}(n)\}$ for distinct $i,j,k$ is equivalent to
the event $\{X_{e_1}(n)=X_{e_2}(n)=X_{e_3}(n)\}$. Notice also that
$P(E_{c,2}(n=2))=1$. So, for $n>2$, there are two or three distinct
edge extrema with probability 1. Hence
$P(E_{c,3}(n))+P(E_{c,2}(n))=1$ for $n > 2$.

By simple integral calculus, we can show that $P(E_{c,2}(n)) \rightarrow 0$ as $n \rightarrow
\infty$, which will imply the desired result.
$\blacksquare$

\subsubsection*{Proof of Theorem \ref{thm:gnNYr=3} }
 Note that $(\Tr)^o
\not= \emptyset$ iff $r<3/2$. Suppose $M \in (\Tr)^o$. Then for any
point $u$ in $R_M(\y_j)$, $\NPE^r(u,M) \subsetneq \TY$, because
there is a tiny strip adjacent to edge $e_j$ not covered by
$\NPE^r(u,M)$, for each $j \in \{1,2,3\}$. Then, $\NPE^r(u,M) \cup
\NPE^r(v,M) \subsetneq \TY$ for all $(u,v) \in R_M(\y_1) \times R_M(\y_2)$.
Pick
$\sup\text{}_{(u,v) \in R_M(\y_1) \times R_M(\y_2)}
\NPE^r(u,M) \cup \NPE^r(v,M)\subsetneq \TY.$
Then
$\TY \setminus \left[\sup\text{}_{(u,v) \in R_M(\y_1) \times R_M(\y_2)}
\NPE^r(u,M) \cup \NPE^r(v,M)\right] $ has positive area.
So
$$\X_n \cap \left[\TY \setminus
\left[\sup\text{}_{(u,v) \in R_M(\y_1) \times R_M(\y_2)} \NPE^r(u,M)
\cup \NPE^r(v,M)\right]\right]\not= \emptyset$$
 with probability 1 for sufficiently large $n$.
(The supremum of a set functional $A(x)$ over a range $B$ is defined
as the set $S:=\sup_{x \in B}A(x)$ such that $S$ is the smallest set
satisfying $A(x)\subseteq S$ for all $x \in B$.) Then at least three
points---one for each vertex region--- are required to dominate
$\X_n$. Hence for sufficiently large $n$, $\g_n(r,M) \ge 3$ with
probability 1, but $\kappa\left( \NPE^r\right)=3$ by Theorem
\ref{thm:kappaNYr=3}. Then $\lim_{n \rightarrow \infty}P\left(
\g_n(r,M)=3 \right)=1$
 for $r<3/2$.
$\blacksquare$

\subsubsection*{Proof of Theorem \ref{thm:MinboundaryTr-g>1}}
\label{sec:MinboundaryTr-g>1}
Let $M=(m_1,m_2) \in \partial (\Tr)$,
say $M \in q_3(r,x)$ (recall that $q_j(r,x)$ are defined such that
$d(\y_j,e_j)=r\cdot d(q_j(r,x),\y_j)$ for $j \in \{1,2,3\}$), then
$m_2=\frac{\sqrt{3}\,(2-r)}{2\,r}$ and $m_1 \in \left[
\frac{3\,(r-1)}{2\,r},\frac{3-r}{2\,r} \right]$. Let $X_{e_j}(n)$ be
one of the closest point(s) to the edge $e_j$; i.e., $X_{e_j}(n)\in
\argmin_{X \in \X_n}d(X,e_j)$ for $j \in \{1,2,3\}$. Note that
$X_{e_j}(n)$ is unique a.s. for each $j$.

Notice that for all $j \in \{1,2,3\}$, $X_{e_j}(n) \notin \NPE^r(X)$
for all $X \in \X_n \cap R_M(\y_j)$ implies that $\g_n(r,M)>1$ with
probability 1. For sufficiently large $n$, $X_{e_j}(n) \notin
\NPE^r(X)$ for all $X \in \X_n \cap R_M(\y_j)$ with probability 1,
for $j \in \{1,2\}$, by the choice of $M$. Hence we consider only
$X_{e_3}(n)$. The asymptotically accurate pdf of $X_{e_3}(n)$ is
$$f_{e_3}\,(x,y)=n\left(\frac{A(S_U(x,y))}{A(\TY)}\right)^{n-1}\frac{1}{A(\TY)},$$
where $S_U(x,y)$ is the unshaded region in Figure \ref{fig:M-in-bound-of-Tr1} (left)
(for a given $X_{e_3}(n)=x_{e_3}=(x,y)$)  whose area is
$A(S_U(x,y))=\sqrt{3}\,\left( 2\,y-\sqrt{3} \right)^2/12$.
Note that $X_{e_3}(n) \notin \NPE^r(X)$ for all
$X \in \X_n \cap R_M(\y_3)$ iff
$\X_n \cap [\G^r_1\left( \X_n,M\right)\cap R_M(\y_3)]=\emptyset$.
Then given $X_{e_3}(n)=(x,y)$,
{\small
$$P\left( \X_n \cap \bigl[ \G^r_1\left( \X_n,M\right)\cap R_M(\y_3)\bigr] =\emptyset \right)=
\left(\frac{A(S_U(x,y))-A\left( \G^r_1\left( \X_n,M\right)\cap R_M(\y_3)\right)}
{A(S_U(x,y))}\right)^{n-1},$$
}
where
$A\left( \G^r_1\left( \X_n,M\right)\cap R_M(\y_3)\right)=
\frac{\sqrt{3}\,y^2}{3\,(r-1)\,r}$
(see Figure \ref{fig:M-in-bound-of-Tr1} (right) where the
shaded region is $\G^r_1\left( \X_n,M\right)\cap R_M(\y_3)$
for a given $X_{e_3}(n)=(x,y)$),
then for sufficiently large $n$
\begin{multline*}
P\left( \X_n \cap \left[ \G^r_1\left( \X_n,M\right)\cap R_M(\y_3)\right] =
\emptyset \right) \approx\\
\int\left(\frac{A(S_U(x,y))-A\left( \G^r_1\left( \X_n,M\right)
\cap R_M(\y_3)\right)}{A(S_U(x,y))}\right)^{n-1}f_{e_3}\,(x,y)\,dy\,dx\\
=\int \frac{n}{A(\TY)}\, \left(\frac{A(S_U(x,y))-
A\left( \G^r_1\left( \X_n,M\right)\cap R_M(\y_3)\right)}{A(\TY)}\right)^{n-1}
\,dy\,dx.
\end{multline*}

Let
{\small
$$
G(x,y)=\frac{A(S_U(x,y))-A\left( \G^r_1\left( \X_n,M\right)\cap R_M(\y_3)\right)}
{A(\TY)} =
\frac{4}{\sqrt{3}}\left(\frac{\sqrt{3}\,\left( 2\,y-\sqrt{3} \right)^2}{12}-
\frac{\sqrt{3}\,y^2}{3\,(r-1)\,r}\right),
$$
}
which is independent on $x$, so we denote it as $G(y)$.
\begin{figure} [ht]
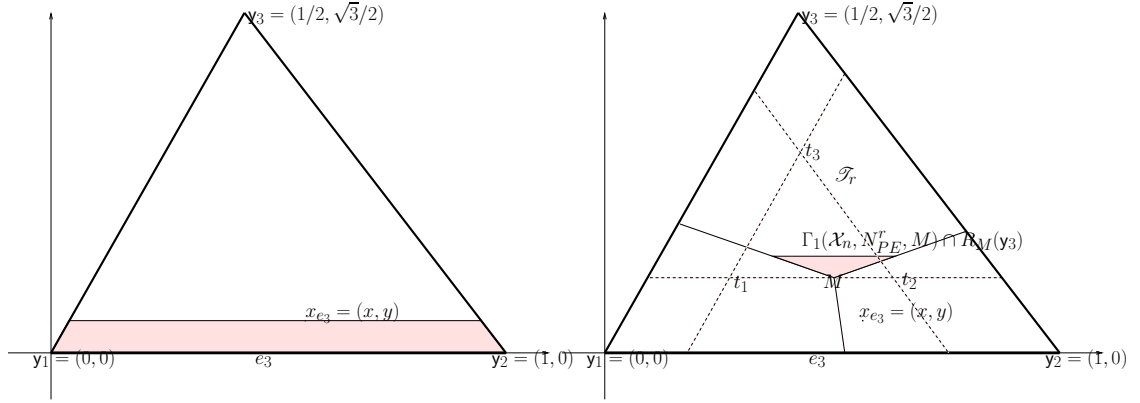

   \centering
    \scalebox{.3}{\input{Xe3pdf.pstex_t}}
    \scalebox{.3}{\input{MinboundaryofTr.pstex_t}}
   \caption{A figure for the description of the pdf of $X_{e_3}(n)$ (left) and
$\G^r_1\left( \X_n,M\right)$ (right) given $X_{e_3}(n)=x_{e_3}=(x,y)$.}
\label{fig:M-in-bound-of-Tr1}
\end{figure}

Let $\ve>0$ be sufficiently small, then for sufficiently large $n$,
\begin{multline*}
P\left( \X_n \cap \left[ \G^r_1\left( \X_n,M\right)\cap R_M(\y_3)\right] =
\emptyset \right) \approx \\
\int_0^{\ve} \int_{y/\sqrt{3}}^{1-y/\sqrt{3}} n\, G(y)^{n-1}\,4/\sqrt{3}\,dy\,dx
= \left( 1-2\,y/\sqrt{3} \right) \int_0^{\ve}n\, G(y)^{n-1}\,4/\sqrt{3}\,dy.
\end{multline*}
 The integrand is critical at $y=0$, since $G(0)=1$ (i.e., when $x_{e_3}\in e_3$).
Furthermore, $G(y)= 1-4\,y/\sqrt{3}+O\left( y^2 \right)$ around $y=0$.
Then letting $y=w/n$, we get
{\small
\begin{eqnarray*}
P\left( \X_n \cap \left[ \G^r_1\left( \X_n,M\right)\cap R_M(\y_3)\right] =
\emptyset \right)
&\approx&
\left( 1-\frac{2\,w}{\sqrt{3}\,n} \right)\,\frac{4}{\sqrt{3}}\,
\int_0^{n\ve} \left( 1-\frac{4\,w}{\sqrt{3}\,n}+ O\left( n^{-2} \right)\right)^{n-1}dw.\\
\text{letting $n \rightarrow \infty$,}
& &
\approx 4/\sqrt{3}\,\int_0^{\infty} \,\exp\left( -4\,w/\sqrt{3} \right)\,dw=1.
\end{eqnarray*}
} Hence $\lim_{n \rightarrow \infty} P\left( \g_n(r,M)>1 \right)=1$.
For $M \in q_j(r,x) \cap \Tr$ with $j \in \{1,2\}$ the result
follows similarly. $\blacksquare$

\subsubsection*{Proof of Theorem \ref{thm:g_nNYr=3-for-r<3/2}}
\label{sec:g_nNYr=3-for-r<3/2} Let $M=(m_1,m_2) \in
\partial(\Tr)\setminus \{t_1(r),\,t_2(r),\,t_3(r)\}$, say $M \in
q_3(r,x)$. Then $m_2=\frac{\sqrt{3}\,(r-1)}{2\,r}$. Without loss of
generality, assume $\frac{1}{2} \le m_1 <  \frac{3-r}{2\,r}$. See
also Figure \ref{fig:M-in-bound-of-Tr2}.

\begin{figure} [ht]
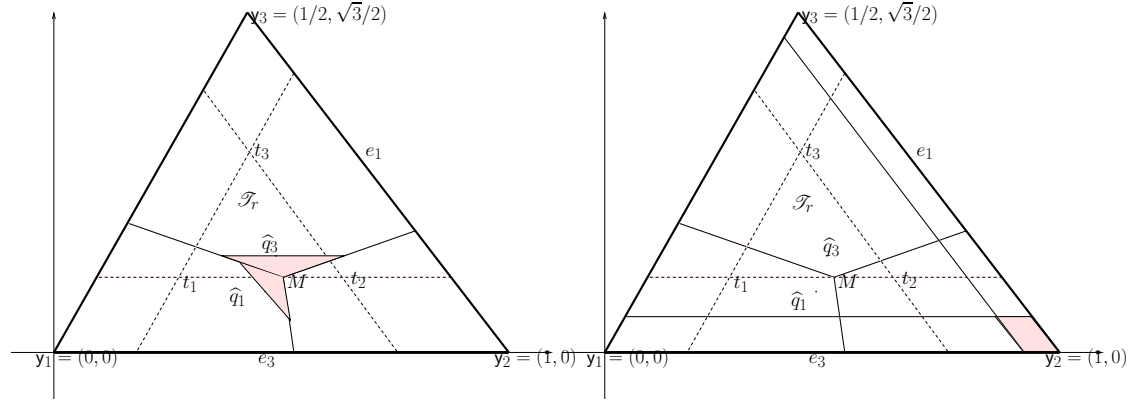

   \centering
   \scalebox{.3}{\input{XM1andXM3pdf.pstex_t}}
   \scalebox{.3}{\input{NYofXM1UXM3.pstex_t}}
   \caption{A figure for the description of the pdf of $\hat{Q}_1(n)$ and
$\hat{Q}_3(n)$ (left) and the unshaded region is
$\NPE^r(\hat{q}_1,M)\cup \NPE^r(\hat{q}_3,M)$ (right).}
\label{fig:M-in-bound-of-Tr2}
\end{figure}

Whenever $\X_n \cap R_M(\y_j)\not= \emptyset$, let
$$\widehat{Q}_j(n)\in \argmin_{X \in \X_n \cap R_M(\y_j)}d\left( X,e_j \right)=\argmax_{X \in \X_n \cap R_M(\y_j)}d(\ell(\y_j,X),\y_j)\text{ for }j \in \{1,2,3\}.$$
Note that at least one of the $\widehat{Q}_j(n)$ uniquely exists w.p. 1 for finite $n$
and as $n \rightarrow \infty$, $\widehat{Q}_j(n)$ are unique w.p. 1.
Then
\begin{multline*}
\g_n(r,M) \le 2 \text{ iff }\X_n \subset \left[ \NPE^r\left( \widehat{Q}_1(n),M \right)  \cup \NPE^r\left( \widehat{Q}_2(n),M \right) \right] \text{ or }\\
\X_n \subset \left[ \NPE^r\left( \widehat{Q}_2(n),M \right)  \cup \NPE^r\left( \widehat{Q}_3(n),M \right) \right] \text{ or }\X_n \subset \left[ \NPE^r\left( \widehat{Q}_1(n),M \right)  \cup \NPE^r\left( \widehat{Q}_3(n),M \right) \right].
\end{multline*}
  Let $E^{i,j}_n$ be the event that
$\X_n \subset \NPE^r\left( \widehat{Q}_i,M \right)  \cup \left[ \NPE^r\left( \widehat{Q}_j(n),M \right) \right]$
for $(i,j) \in \bigl\{ (1,2),(1,3),(2,3) \bigr\}$.
Then
\begin{multline*}
P\left( \g_n(r,M) \le 2 \right)= P\left( E^{1,2}_n \right)+P\left(
E^{2,3}_n \right)+P\left( E^{1,3}_n \right)-
P\left( E^{1,2}_n \cap E^{2,3}_n \right)-P\left( E^{1,2}_n \cap E^{1,3}_n \right)\\
-P\left( E^{1,3}_n \cap E^{2,3}_n \right)+
P\left( E^{1,2}_n \cap E^{2,3}_n\cap E^{1,3}_n \right).
\end{multline*}
But note that $P\left( E^{1,2}_n \right)\rightarrow 0$ as
$n \rightarrow \infty$ by the choice of $M$ since
$$\sup\text{}_{\substack{u \in R_M(\y_1)\\v \in R_M(\y_2)}}\NPE^r(u,M) \cup \NPE^r(v,M) \subsetneq \TY,$$
 and
$$P\left( \X_n \cap \TY \setminus
\left[ \sup\text{}_{ \substack{u\in R_M(\y_1)\\v \in R_M(\y_2)}}\NPE^r(u,M) \cup \NPE^r(v,M) \right]
\not= \emptyset \right)
\rightarrow 1 \text{ as } n \rightarrow \infty.$$
Then,
$$P\left( E^{1,2}_n \right)-P\left( E^{1,2}_n \cap E^{2,3}_n \right)-P\left( E^{1,2}_n \cap E^{1,3}_n \right)+P\left( E^{1,2}_n \cap E^{2,3}_n\cap E^{1,3}_n \right) \le 4\,P\left( E^{1,2}_n \right) \rightarrow 0 \text{ as } n \rightarrow \infty.$$
Therefore,
$$\lim_{n \rightarrow \infty}P\left( \g_n(r,M) \le 2 \right)=\lim_{n \rightarrow \infty} \left( P\left( E^{2,3}_n \right)+P\left( E^{1,3}_n \right) \right).$$

Furthermore, observe that
$P\left( E^{1,3}_n \right) \ge P\left( E^{2,3}_n \right)$ by the choice of $M$.
Then we first find $\lim_{n \rightarrow \infty} P\left( E^{1,3}_n \right)$.
Given a realization of $\X_n$ with $\widehat{Q}_1(n)=\widehat{q}_1=(x_1,y_1)$ and
$\widehat{Q}_3(n)=\widehat{q}_3=(x_3,y_3)$, the remaining $n-2$ points should fall,
for example, in the undshaded region in
Figure \ref{fig:M-in-bound-of-Tr2} (left).
Then the asymptotically accurate joint pdf of $\widehat{Q}_1(n),\widehat{Q}_3(n)$ is
$$f_{13}\,\bigl( \vec{\zeta} \bigr)=
\frac{n\,(n-1)}{A(\TY)^2}\left(\frac{A(\TY)-A(S_R\bigl( \vec{\zeta} \bigr)
\bigl( \vec{\zeta} \bigr))}{A(\TY)}\right)^{n-2}$$
 where $\vec{\zeta}=(x_1,y_1,x_3,y_3)$, $S_R\bigl( \vec{\zeta} \bigr)$ is the shaded region
in Figure \ref{fig:M-in-bound-of-Tr2} (left) whose area is
$A(S_R\bigl( \vec{\zeta} \bigr))=
\frac{\sqrt{3}\,\left( 2\,r\,y_3-\sqrt{3}\,(r-1) \right)^2}{12\,r\,(r-1)}+
\frac{\sqrt{3}\left[ 2\,\sqrt{3}\,r\,y_1-3\,(r-1)+6\,r\,(x_1-m_1) \right]^2}
{72\,r\,(1-r\,(2\,m_1-1))}$.

Given $\widehat{Q}_j(n)=\widehat{q}_j=(x_j,y_j)$ for $j \in \{1,3\}$,
$$P\left( E^{1,3}_n \right)=
\left(\frac{A\left( \NPE^r\left( \widehat{q}_1,M \right)  \cup \NPE^r
\left( \widehat{q}_3,M \right)  \right)-
A(S_R\bigl( \vec{\zeta} \bigr))}
{A(\TY)-A(S_R\bigl( \vec{\zeta} \bigr))}\right)^{n-2}$$
 then for sufficiently large $n$
\begin{eqnarray*}
P\left( E^{1,3}_n \right)&\approx&\int \left(\frac{A\left( \NPE^r\left( \widehat{q}_1,M \right)  \cup \NPE^r\left( \widehat{q}_3,M \right)  \right)-A(S_R\bigl( \vec{\zeta} \bigr))}{A(\TY)-A(S_R\bigl( \vec{\zeta} \bigr))}\right)^{n-2} f_{13}\,\bigl( \vec{\zeta} \bigr)\,d\vec{\zeta},\\
 & = &\int \frac{n\,(n-1)}{A(\TY)^2}\left(\frac{A\left( \NPE^r\left( \widehat{q}_1,M \right)  \cup \NPE^r\left( \widehat{q}_3,M \right)  \right)-A(S_R\bigl( \vec{\zeta} \bigr))}{A(\TY)}\right)^{n-2}\,d\vec{\zeta}
\end{eqnarray*}
where
$$A\left( \NPE^r\left( \widehat{q}_1,M \right)  \cup \NPE^r\left( \widehat{q}_3,M \right)  \right)=\frac{\sqrt{3}}{4}-\left(\frac{\left( \sqrt{3}\,r\,y_1+3\,r\,x_1-3 \right)\,\left( \sqrt{3}\,(r-1)-2\,r\,y_3 \right)}{6}\right).$$
See Figure \ref{fig:M-in-bound-of-Tr2} (right) for
$\NPE^r(\widehat{q}_1,M) \cup \NPE^r(\widehat{q}_3,M)$.
Let
$$G\bigl( \vec{\zeta} \bigr)=\frac{A\left( \NPE^r\left( \widehat{q}_1,M \right)  \cup \NPE^r\left( \widehat{q}_3,M \right)  \right)-A(S_R\bigl( \vec{\zeta} \bigr))}{A(\TY)}.$$
Note that the integral is critical at $x_1=x_3=m_1$ and $y_1=y_3=m_2$,
since $G\bigl( \vec{\zeta} \bigr)=1$.
Since $\NPE^r(x,M_C)$ depends on the distance $d(x,e_j)$ for $x \in R_M(\y_j)$,
we make the change of variables $(x_1,y_1) \rightarrow (d(M,e_1)+z_1,y_1)$
where $d(M,e_1)=\frac{\sqrt{3}\,(r+1-2\,r\,m_1)}{4\,r}$ and
$(x_3,y_3) \rightarrow (x_3,m_2+z_3)$ then $G\bigl( \vec{\zeta} \bigr)$
depends only on $z_1,z_3$, we denote it $G(z_1,z_3)$ which is
$$G(z_1,z_3)=1-\frac{8\,r\,z_1^2}{3\,(1+r\,(1-2\,m_1))}-\frac{4\,r\,z_3^2}{3\,(r-1)}-\frac{2\,r\,z_3\,\left( \sqrt{3}\,(3-r) \right)+r\,\left( 4\,z_1-2\,\sqrt{3}\,m_1 \right)}{3}.$$
  The new integrand is $\frac{n\,(n-1)}{A(\TY)^2}G(z_1,z_3)^{n-2}$.
Integrating with respect to $x_3$ and $y_1$ yields $\frac{2\,\sqrt{3}\,z_3\,r}{3\,(r-1)}$
and $\frac{4\,\sqrt{3}\,r\,z_1}{3\,(2\,r\,m_1-r-1)}$, respectively.
Hence for sufficiently large $n$
$$
P\left( E^{1,3}_n \right) \approx \int_0^{\ve}\int_0^{\ve}\,\frac{n\,(n-1)}{A(\TY)^2}\,\left(\frac{2\,\sqrt{3}\,z_3\,r}{3\,(r-1)}\right)\,\left(\frac{4\,\sqrt{3}\,r\,z_1}{3\,(2\,r\,m_1-r-1)}\right) G(z_1,z_3)^{n-2}dz_1\,dz_3.
$$
 Note that the new integral is critical when $z_1=z_3=0$,
so we make the change of variables $z_1=w_1/\sqrt{n}$ and
$z_3=w_3/n$ then $G(z_1,z_3)$ becomes
$$G(w_1,w_3)=1+\frac{1}{n}\,\left(\frac{2\,\sqrt{3}\,r\,(r-3+2\,r\,m_1)}{3}w_3+\frac{8\,r}{3\,(r+1-2\,r\,m_1)}\,w_1^2\right)+O\left( n^{-3/2} \right),$$
so for sufficiently large $n$
\begin{align*}
&P\left( E^{1,3}_n \right) \approx \int_0^{\sqrt{n}\,\ve}\int_0^{n\,\ve}\,\frac{(n-1)}{n^3}\,\frac{16}{3}\,\left(\frac{2\,\sqrt{3}\,r}{3\,(r-1)}\right)\,\left(\frac{4\,\sqrt{3}\,r}{3\,(2\,r\,m_1-r-1)}\right)\,(-4\,m_1+2+\sqrt{2}) w_1\,w_3\\
& \,\Biggl[1-\frac{1}{n}\Biggl(\frac{2\,\sqrt{3}\,r\,(r-3+2\,r\,m_1)}{3}w_3+\frac{8\,r}{3\,(r+1-2\,r\,m_1)}\,w_1^2\Biggr)+O\left( n^{-3/2} \right)\Biggr]^{n-2}dw_3w_1,\\
&\approx O\left( n^{-1} \right) \int_0^{\infty}\int_0^{\infty}\,w_1\,w_3\,
\exp\left(-\frac{2\,\sqrt{3}\,r\,(r-3+2\,r\,m_1)w_3}{3}-\frac{8\,r\,w_1^2}
{3\,(r+1-2\,r\,m_1)}\right)\,dw_3w_1 = O\left( n^{-1} \right)
\end{align*}
 since
$\int_0^{\infty}\int_0^{\infty}\,w_1\,w_3\,\exp\left(-\frac{2\,\sqrt{3}\,r\,(r-3+2\,r\,m_1)}{3}w_3-\frac{8\,r}{3\,(r+1-2\,r\,m_1)}\,w_1^2\right)\,dw_3w_1=\frac{3}{8\,r\,(3-r\,(2\,m_1+1))},$
 which is a finite constant.
Then $P\left( E^{1,3}_n \right) \rightarrow 0$ as $n \rightarrow
\infty$, which also implies $P\left( E^{2,3}_n \right) \rightarrow
0$ as $n \rightarrow \infty$. Then $P\left( \g_n(r,M) \le 2 \right)
\rightarrow 0$. Hence the desired result follows. $\blacksquare$

\subsubsection*{Proof of Theorem \ref{thm:g_nNYr=2-for-r-vertex-of-Tr}}
\label{sec:g_nNYr=2-for-r-vertex-of-Tr} Let $M=(m_1,m_2) \in
\{t_1(r),\,t_2(r),\,t_3(r)\}$. Without loss of generality, assume
$M=t_2(r)$ then $m_1=\frac{2-r+c_1\,(r-1)}{r}$ and
$m_2=\frac{c_2\,(r-1)}{r}$. See Figure \ref{fig:M-in-bound-of-Tr3}.

\begin{figure} [ht]
   \centering
   \scalebox{.3}{\input{XM1andXM3pdf_Mt2.pstex_t}}
   \scalebox{.3}{\input{NYofXM1UXM3_Mt2.pstex_t}}
\caption{
A figure for the description of the pdf of $\hat{Q}_1(n)$ and $\hat{Q}_3(n)$ (left)
and the unshaded region is $\NPE^r(\hat{q}_1,M) \cup \NPE^r(\hat{q}_3,M)$ (right)
given $\hat{Q}_j(n)=\hat{q}_j$ for $j \in \{1,3\}$.}
\label{fig:M-in-bound-of-Tr3}
\end{figure}

Let $\widehat{Q}_j(n)$ and the events $E^{i,j}_n$ be defined as in the proof of
Theorem \ref{thm:g_nNYr=3-for-r<3/2} for
$(i,j) \in \bigl\{ (1,2),(1,3),(2,3) \bigr\}$.
Then as in the proof of Theorem \ref{thm:g_nNYr=3-for-r<3/2},
\begin{multline*}
P\left( \g_n(r,M) \le 2 \right)= P\left( E^{1,2}_n \right)+ P\left(
E^{2,3}_n \right)+ P\left( E^{1,3}_n \right)-
P\left( E^{1,2}_n \cap E^{2,3}_n \right)-\\
P\left( E^{1,2}_n \cap E^{1,3}_n \right)-
P\left( E^{1,3}_n \cap E^{2,3}_n \right)+
P\left( E^{1,2}_n \cap E^{2,3}_n\cap E^{1,3}_n \right).
\end{multline*}

Observe that the choice of $M$ implies that
$P\left( E^{1,3}_n \right) \ge P\left( E^{2,3}_n \right)$ and
by symmetry (in $T_e$) $P\left( E^{1,2}_n \right)=P\left( E^{2,3}_n \right)$.
So first we find $P\left( E^{1,3}_n \right)$.
As in the proof of Theorem \ref{thm:g_nNYr=3-for-r<3/2}
asymptotically accurate joint pdf of $\widehat{Q}_1(n),\widehat{Q}_3(n)$ is
$$f_{13}\,\bigl( \vec{\zeta} \bigr)=\frac{n\,(n-1)}{A(\TY)^2}\left(\frac{A(\TY)-A(S_R\bigl( \vec{\zeta} \bigr))}{A(\TY)}\right)^{n-2}$$
 where $\vec{\zeta}=(x_1,\y_1,x_3,\y_3)$ and $S_R\bigl( \vec{\zeta} \bigr)$ is
the shaded region in Figure \ref{fig:M-in-bound-of-Tr3} (left) whose area is
$$A(S_R\bigl( \vec{\zeta} \bigr))=\frac{\sqrt{3}\,\left( 2\,r\,y_3-\sqrt{3}\,(r-1)^2 \right)}{12\,(r-1)\,r}+\frac{\sqrt{3}\,\left( \sqrt{3}\,r\,y_1+3\,x_1\,r-3 \right)^2}{36\,(r-1)\,r}.$$

Given $\widehat{Q}_j(n)=\widehat{q}_j=(x_j,y_j)$  for $j \in \{1,3\}$,
$$P\left( E^{1,3}_n \right)=\left(\frac{A\left( \NPE^r\left( \widehat{q}_1,M \right)  \cup \NPE^r\left( \widehat{q}_3,M \right)  \right)-A(S_R\bigl( \vec{\zeta} \bigr))}{A(\TY)-A(S_R\bigl( \vec{\zeta} \bigr))}\right)^{n-2},$$
then for sufficiently large $n$
\begin{eqnarray*}
P\left( E^{1,3}_n \right)&\approx&\int \left(\frac{A\left( \NPE^r\left( \widehat{Q}=\widehat{q}_1,M \right)  \cup \NPE^r\left( \widehat{q}_3,M \right)  \right)-A(S_R\bigl( \vec{\zeta} \bigr))}{A(\TY)-A(S_R\bigl( \vec{\zeta} \bigr))}\right)^{n-2} f_{13}\,\bigl( \vec{\zeta} \bigr)\,d\vec{\zeta},\\
&=& \int \frac{n\,(n-1)}{A(\TY)^2}\left(\frac{A\left( \NPE^r\left( \widehat{q}_1,M \right)  \cup \NPE^r\left( \widehat{q}_3,M \right)  \right)-A(S_R\bigl( \vec{\zeta} \bigr))}{A(\TY)}\right)^{n-2}\,d\vec{\zeta}
\end{eqnarray*}
 where
$$A\left( \NPE^r\left( \widehat{q}_1,M \right)  \cup \NPE^r\left( \widehat{q}_3,M \right)  \right)=\frac{\sqrt{3}}{4}-\frac{\left( 2\,r\,y_3-\sqrt{3}\,(r-1) \right)\,\left( 3-\sqrt{3}\,r\,y_1-3\,r\,x_1 \right)}{6}.$$
See Figure \ref{fig:M-in-bound-of-Tr3} (right) for
$\NPE^r(\widehat{q}_1,M) \cup \NPE^r(\widehat{q}_3,M)$.
 Let
$$G\bigl( \vec{\zeta} \bigr)=\frac{A\left( \NPE^r\left( \widehat{q}_1,M \right)  \cup \NPE^r\left( \widehat{q}_3,M \right)  \right)-A(S_R\bigl( \vec{\zeta} \bigr))}{A(\TY)}.$$
Note that the integral is critical when $x_1=x_3=m_1$ and $y_1=y_3=m_2$,
since $G\bigl( \vec{\zeta} \bigr)=1$.

As in the proof of Theorem \ref{thm:g_nNYr=3-for-r<3/2},
we make the change of variables $(x_1,y_1)\rightarrow (d(M,e_1)+z_1,y_1)$
where $d(M,e_1)=\frac{\sqrt{3}\,(r-1)}{2\,r}$ and
$(x_3,y_3)\rightarrow (x_3,m_2+z_3)$.
Then $G\bigl( \vec{\zeta} \bigr)$ becomes
$$G(z_1,z_3)=1-\frac{4\,r}{3\,(r-1)}\,z_1^2-\frac{4\,r}{3\,(r-1)}\,z_3^2-\frac{8\,r^2}{3}\,\,z_1\,z_3.$$
The new integral is
$$\int \frac{n\,(n-1)}{A(\TY)^2}G(z_1,z_3)^{n-2}dx_3dy_1dz_3dz_1.$$
Note that $G(z_1,z_3)$ is independent of $y_1,x_3$,
so integrating with respect to $x_3$ and $y_1$ yields
$\frac{2\,\sqrt{3}\,r\,z_1}{3\,(r-1)}$ and
$\frac{2\,\sqrt{3}\,r\,z_3}{3\,(r-1)}$, respectively.
The new integral is critical at $z_1=z_3=0$.
Hence, for sufficiently large $n$ and sufficiently small $\ve>0$,
the integral becomes,
$$P\left( E^{1,3}_n \right)\approx \int_0^{\ve}\int_0^{\ve}\,\frac{n\,(n-1)}{A(\TY)^2}\,\left(\frac{12\,r^2}{9\,(r-1)^2}\right)\,z_1\,z_3\,G(z_1,z_3)^{n-2}\,\,dz_1dz_3.$$
 Since the new integral is critical when $z_1=z_2=0$,
we make the change of variables $z_j=w_j/\sqrt{n}$ for $j \in \{1,3\}$;
then $G(z_1,z_3)$ becomes
$$G(w_1,w_3)=1-\frac{4\,r}{3\,n\,(r-1)}\left( w_1^2+w_3^2+2\,r\,(r-1)\,w_1\,w_3)\right),$$
so
\begin{eqnarray*}
\lefteqn{p_r:=P\left( E^{1,3}_n \right) \approx \int_0^{\sqrt{n}\,\ve}\int_0^{\sqrt{n}\,\ve}\,\frac{(n-1)}{n}\,\frac{16}{3}\,\left(2\,\left(\frac{12\,r^2}{9\,(r-1)^2}\right)w_1\,w_3\right)}\\
& &\left[1-\frac{4\,r}{3\,n\,(r-1)}\,(w_1^2+w_3^2+2\,r\,(r-1)\,w_1\,w_3))\right]^{n-2}dw_3w_1,
\;\; \text{letting $n \rightarrow \infty$,}\\
&\approx & \int_0^{\infty}\int_0^{\infty}\,\frac {64}{9}\,\left(\frac{r}{r-1}\right)^2\,w_1\,w_3\,\,\exp\left(\frac{4\,r}{3\,(r-1)}\,(w_1^2+w_3^2+2\,r\,(r-1)\,w_1\,w_3)\right)\,dw_3w_1
\end{eqnarray*}
 which is not analytically integrable, but $p_r$
can be obtained by numerical integration, e.g.,
$p_{r=\sqrt{2}}\approx .4826$ and $p_{r=5/4}\approx .6514$.

Next, we find $\lim_{n \rightarrow \infty}P\left( E^{2,3}_n \right)$.
The asymptotically accurate joint pdf of $\widehat{Q}_2(n),\widehat{Q}_3(n)$ is
$$f_{23}\,\bigl( \vec{\zeta} \bigr)=
\frac{n\,(n-1)}{A(\TY)^2}\left(\frac{A(\TY)-A(S_R^2\bigl( \vec{\zeta} \bigr))}
{A(\TY)}\right)^{n-2}$$
where $\vec{\zeta}=(x_2,y_2,x_3,y_3)$ and $S_R^2\bigl( \vec{\zeta} \bigr)$
is the shaded
region in Figure \ref{fig:M-in-bound-of-Tr4} (left) whose area is
$$A\left( S_R^2\bigl( \vec{\zeta} \bigr) \right)  =\frac{\sqrt{3}\,\left( 2\,r\,\y_3+\sqrt{3}\,(1-r) \right)}{12\,r\,(r-1)}+\frac{\sqrt{3}\,\left( \sqrt{3}\,r\,\y_2-3\,r\,x_2-3\,r+6 \right)}{36\,(2-r)\,r}.$$

\begin{figure} [ht]
   \centering
   \scalebox{.3}{\input{XM2andXM3pdf.pstex_t}}
   \scalebox{.3}{\input{NYofXM2UXM3.pstex_t}}
   \caption{A figure for the description of the pdf of $\hat{Q}_2(n)$ and $\hat{Q}_3(n)$ (left) and unshaded region is $\NPE^r(\hat{q}_2) \cup \NPE^r(\hat{q}_3)$ (right)
given $\hat{Q}_j(n)=\hat{q}_j$ for $j\in\{2,3\}$.}
\label{fig:M-in-bound-of-Tr4}
\end{figure}

As before,
\begin{eqnarray*}
P\left( E^{2,3}_n \right)&=&\int \left(\frac{A\left( \NPE^r\left( \widehat{q}_2,M \right)
\cup \NPE^r\left( \widehat{q}_3,M \right)  \right)-A(S_R\bigl( \vec{\zeta} \bigr))}
{A(\TY)-A(S_R^2\bigl( \vec{\zeta} \bigr))}\right)^{n-2} f_{23}\,\bigl( \vec{\zeta} \bigr)\,d\vec{\zeta}\\
 & = & \int \frac{n\,(n-1)}{A(\TY)^2}\left(\frac{A\left( \NPE^r\left( \widehat{q}_2,M \right)
 \cup \NPE^r\left( \widehat{q}_3,M \right)  \right)-A(S_R\bigl( \vec{\zeta} \bigr))}
 {A(\TY)}\right)^{n-2}\,d\vec{\zeta},
\end{eqnarray*}
 where
$A\left( \NPE^r\left( \widehat{q}_2,M \right)  \cup \NPE^r\left( \widehat{q}_3,M \right) \right)
=\frac{\sqrt{3}}{4}-\frac{\left( 2\,r\,y_3-\sqrt{3}\,(r-1) \right)\,
\left( 3-\sqrt{3}\,r\,y_2+3\,r\,x_2-3\,r \right)}{6}.$

See Figure \ref{fig:M-in-bound-of-Tr4} (right) for
$\NPE^r(\widehat{q}_2) \cup \NPE^r(\widehat{q}_3,M)$.
Let
$$G\bigl( \vec{\zeta} \bigr)=\frac{A\left( \NPE^r\left( \widehat{q}_2,M \right)  \cup \NPE^r\left( \widehat{q}_3,M \right)  \right)-A(S_R\bigl( \vec{\zeta} \bigr))}{A(\TY)}.$$
Note that the integral is critical when $x_2=x_3=m_1$ and $y_2=y_3=m_2$,
since $G\bigl( \vec{\zeta} \bigr)=1$.

We make the change of variables $(x_3,y_3) \rightarrow (x_3,m_2+z_3)$
and $(x_2,y_2)\rightarrow (d(M,e_2)+z_2,y_2)$ where
$d(M,e_2)=\frac{\sqrt{3}\,(2-r)}{2\,r}$.
Then $G\bigl( \vec{\zeta} \bigr)$ becomes
$$G(z_2,z_3)=1-\frac{4\,r\,z_2^2}{3\,(2-r)}-\frac{4\,r\,z_2^2}{3\,(r-2)}-\frac{4\,\sqrt{3}\,r\,z_3\,(3-2\,r)}{3}-\frac{8\,r^2\,z_2\,z_3}{3}.$$
The new integral is
$$\int \frac{n\,(n-1)}{A(\TY)^2}G(z_2,z_3)^{n-2}dx_3dy_2dz_3dz_2.$$
The integrand is independent of $x_3$ and $y_2$,
so integrating with respect to $x_3$ and $y_2$ yields
$\frac{2\,\sqrt{3}\,r\,z_3}{3\,(r-1)}$ and
$\frac{2\,\sqrt{3}\,r\,z_2}{3\,(2-r)}$, respectively.
Hence, for sufficiently large $n$
$$P\left( E^{2,3}_n \right) \approx\int_0^{\ve}\int_0^{\ve}\,\frac{n\,(n-1)}{A(\TY)^2}\,\left(\frac{4\,r^2}{3\,(r-1)\,(2-r)}\right)z_3\,z_2\,G(z_2,z_3)^{n-2}dz_2dz_3.$$
Note that the new integral is critical when $z_2=z_3=0$,
so we make the change of variables $z_2=w_2/\sqrt{n}$
and $z_3=w_3/n$ then $G(z_2,z_3)$ becomes
$$G(w_2,w_3)=1-\frac{1}{n}\left[ \frac{4\,r\,w_2^2}{3\,(2-r)}-\frac{4\,\sqrt{3}\,r\,w_3\,(3-2\,r)}{3} \right]+O\left( n^{-\frac{3}{2}} \right),$$
so for sufficiently large $n$
\begin{align*}
P\left( E^{2,3}_n \right)&
\approx \int_0^{\sqrt{n}\,\ve}\int_0^{n\,\ve}\,\frac{(n-1)}{n^2}\,\frac{64\,r^2}
{9\,(r-1)\,(2-r)}\,w_2\,w_3\\
& \left[1-\frac{1}{n}\left(\frac{4\,r\,w_2^2}{3\,(2-r)}-
\frac{4\,\sqrt{3}\,r\,w_3\,(3-2\,r)}{3}\right)+
O\left( n^{-\frac{3}{2}} \right)\right]^{n-2}dw_3w_2,\\
&\approx O\left( n^{-1} \right)\,\int_0^{\infty}\int_0^{\infty}
\,w_2\,w_3\,\,\exp\left(-\frac{4\,r\,w_2^2}{3\,(2-r)}-
\frac{4\,\sqrt{3}\,r\,u_3\,(3-2\,r)}{3}\right)\,dw_3w_2 =
O\left( n^{-1} \right)
\end{align*}
 since
$$\int_0^{\infty}\int_0^{\infty}\, w_2\,w_3\,\,\exp\left(-\frac{4\,r\,w_2^2}{3\,(2-r)}-\frac{4\,\sqrt{3}\,r\,u_3\,(3-2\,r)}{3}\right)\,dw_3w_2=\frac{27\,(2-r)}{384\,r^3\,(3-2\,r)^2}$$
which is a finite constant.

Thus we have shown that $P\left( E^{2,3}_n \right) \rightarrow 0$
as $n \rightarrow \infty$, which implies that as $n \rightarrow \infty$,
\begin{multline*}
P\left( E^{2,3}_n \right)+P\left( E^{1,2}_n \right)-P\left( E^{1,2}_n \cap E^{2,3}_n \right)-P\left( E^{1,2}_n \cap E^{1,3}_n \right)\\
-P\left( E^{1,3}_n \cap E^{2,3}_n \right)+P\left( E^{1,2}_n \cap E^{2,3}_n\cap E^{1,3}_n \right) \le 5\,P\left( E^{2,3}_n \right) \rightarrow 0.
\end{multline*}

Hence $\lim_{n \rightarrow \infty}P\left( \g_n(r,M) \le 2
\right)=\lim_{n \rightarrow \infty}P\left( E^{1,3}_n \right)$ and
$\lim_{n \rightarrow \infty}P(\g_n(r,M) >1)=1$ together imply that
$$\lim_{n \rightarrow \infty}P(\g_n(r,M) = 2)=p_r.\; \blacksquare$$
}

\end{document}